\documentclass[11pt]{amsart}

\usepackage{amscd}
\usepackage{amssymb}
\usepackage{amsmath}
\usepackage{amsfonts}

\usepackage[all]{xy}

\newtheorem{proposition}{Proposition}[section]

\newtheorem{theorem}[proposition]{Theorem}
\newtheorem{cor}[proposition]{Corollary}

\newcounter{THNO}[section]
\newcounter{SNO}[section]
\newcounter{multieq}[equation]

\newcounter{tmp}

\def\db#1{ {\mathbf D}^b({#1})}

\def\h#1,#2{{\rm Hom}({#1}\:,\; {#2})}
\def\H#1,#2,#3,#4{{\rm Hom}^{#1}_{#2}({#3}\:,\; {#4})}
\def\E#1,#2,#3,#4{{\rm Ext}^{#1}_{#2}({#3}\:,\; {#4})}
\def\lto{\longrightarrow}

\def\da{\big\downarrow}

\def\ot{\otimes}

\def\ts{\times}
\def\op{\oplus}

\def\wh#1{\widehat{#1}}

\def\C{{\mathbb{C}}}

\def\llongrightarrow{{\;\relbar\joinrel\relbar\joinrel\relbar\joinrel
\rightarrow\;}}

\def\A{{\mathcal{A}}}
\def\O{{\mathcal{O}}}
\def\E{{\mathcal{E}}}
\def\F{{\mathcal{F}}}

\def\H{{\mathcal{H}}}

\def\J{{\mathcal{J}}}

\def\End{{\rm E}{\rm n}{\rm d}}
\def\Hom{{\rm H}{\rm o}{\rm m}}

\def\Pic{{\rm P}{\rm i}{\rm c}}

\def\h{{\hbar}}

\def\Endo{{{\mathcal E}nd}}
\def\X{{\mathcal{X}}}

\def\be{\begin{equation}}
\def\ee{\end{equation}}

\def\coh{{\rm coh}}
\def\bdot{\begin{picture}(4,4)\put(2,3){\circle*{1.5}}\end{picture}}
\def\bH{{\mathbf H}}
\def\bC{{\mathbf C}}
\def\bD{{\mathbf D}}

\newcommand{\CC}{{\mathbb C}}
\newcommand{\RR}{{\mathbb R}}
\newcommand{\ZZ}{{\mathbb Z}}
\newcommand{\PP}{{\mathbb P}}

\newcommand{\QQ}{{\mathbb Q}}

\newcommand{\LY}{{{\mathcal L}Y}}
\newcommand{\NY}{{{\mathcal N}Y}}
\newcommand{\fH}{{\mathfrak H}}

\newcommand{\ra}{\rightarrow}

\newcommand{\hA}{{\hat A}}

\newcommand{\id}{\mathrm{id}}

\newcommand{\Ker}{\mathop{\mathrm{Ker}}\nolimits}
\renewcommand{\Im}{\mathop{\mathrm{Im}}\nolimits}

\newcommand{\Gi}{{\left(G^{-1}\right)}}
\newcommand{\vac}{{|vac\rangle}}
\newcommand{\ds}{\displaystyle}

\newcommand{\bz}{{\bar{z}}}

\newcommand{\bT}{{\bar{T}}}
\newcommand{\ba}{{\bar{a}}}

\newcommand{\btheta}{{\bar{\theta}}}
\newcommand{\al}{{\alpha}}
\newcommand{\bal}{{\bar{\alpha}}}
\newcommand{\tpsi}{{\bar{\psi}}}
\newcommand{\hw}{{{W}}}
\newcommand{\hm}{{{M}}}
\newcommand{\bpartial}{{\bar{\partial}}}
\newcommand{\hp}{{{P}}}
\newcommand{\htp}{{{\bar{P}}}}

\newcommand{\grp}{{\mathbb C}\ [\Gamma\op\Gamma^*]}

\newcommand{\cI}{{\mathcal I}}
\newcommand{\cJ}{{\mathcal J}}

\newcommand{\cB}{{\mathcal B}}
\newcommand{\cH}{{\mathcal H}}

\newcommand{\cV}{{\mathcal V}}
\newcommand{\cF}{{\mathcal F}}

\newcommand{\B}{B}

\newcommand{\m}{{\mathfrak M}}

\newcommand{\bpsi}{{\bar\psi}}

\newcommand{\bL}{{\bar L}}
\newcommand{\bQ}{{\bar Q}}
\newcommand{\bJ}{{\bar J}}

\newcommand{\SCFT}{{\rm SCFT}}
\newcommand{\ch}{{\rm ch}}

\title[]{Lectures on Mirror Symmetry, Derived Categories, and D-branes}

\author[]{Anton Kapustin}

\address{California Institute of Technology, Pasadena, CA 91125, USA}
\email{kapustin@theory.caltech.edu}

\author[]{Dmitri Orlov}

\address{Algebra Section, Steklov Mathematical Institute,
Russian Academy of Sciences, 8 Gubkin str., Moscow 119996, Russia}
\email{orlov@mi.ras.ru}

\thanks{The first author was supported in part by the DOE grant DE-FG03-92-ER40701.
The second author was supported in part by the grant of the
President of RF for young scientists ³D-2731.2004.1,  by Civic
Research Development Foundation (CRDF Award No RM1-2405-MO-02),
and by the Russian Science Support Foundation.}

\begin{document}

\begin{abstract}
This paper is  an introduction to Homological Mirror Symmetry,
derived categories, and topological D-branes aimed mainly at
a mathematical audience. In the paper we explain the physicists'
viewpoint of the Mirror Phenomenon, its relation to derived
categories, and the reason why it is necessary to enlarge the
Fukaya category with coisotropic A-branes; we  discuss how to
extend the definition of Floer homology to such objects and
describe mirror symmetry for flat tori. The paper consists of four
lectures which were given at the Institute for Pure and Applied Mathematics
(Los Angeles), March 2003, as part of a program on Symplectic
Geometry and Physics.
\end{abstract}

\maketitle

\section{Mirror Symmetry From a Physical Viewpoint}

The goal of the first lecture is to explain the physicists'
viewpoint of the Mirror Phenomenon and its interpretation in
mathematical terms proposed by Maxim Kontsevich in his 1994 talk
at the International Congress of Mathematicians~\cite{Konts}.
Another approach to Mirror Symmetry was proposed by A.~Strominger,
S-T.~Yau, and E.~Zaslow~\cite{SYZ}, but we will not discuss it
here.

From the physical point of view, Mirror Symmetry is a relation on
the set of 2d conformal field theories with $N=2$ supersymmetry. A
2d conformal field theory is a rather complicated algebraic object
whose definition will be sketched in a moment. Thus Mirror
Symmetry originates in the realm of algebra. Geometry will appear
later, when we specialize to a particular class of $N=2$
superconformal field theories related to Calabi-Yau manifolds.

Let us start with 2d conformal field theory. The data needed to
specify a 2d CFT consist of an infinite-dimensional vector space
$V$ (the space of states), three special elements in $V$ (the
vacuum vector $\vac,$ and two more elements $L$ and $\bL$), and a
linear map $Y$ from $V$ to the space of ``formal fractional power
series in $z,\bz$ with coefficients in $\End(V)$'' ($Y$ is called
the state-operator correspondence). The precise definition of what
a ``formal fractional power series'' means can be found in
~\cite{KO}; to keep things simple, one can pretend that $Y$ takes
values in the space of Laurent series in $z,\bz$ with coefficients
in $\End(V),$ although such a definition is not sufficient for
applications to Mirror Symmetry. These data must satisfy a number
of axioms whose precise form can be found in ~\cite{KO}. Roughly
speaking, they are

\begin{itemize}
\item[(i)] $Y(\vac)=\id_V.$ \item[(ii)] $Y(L)=\sum_{n\in\ZZ}
\frac{L_n}{z^{n+2}},\quad Y(\bL)=\sum_{n\in\ZZ}
\frac{\bL_n}{\bz^{n+2}}$ for some $L_n,\bL_n\in \End(V).$
\item[(iii)] Both $L_n$ and $\bL_n$ satisfy the commutation
relations of the Virasoro algebra, and all $L_n$ commute with all
$\bL_m.$ \item[(iv)] $[L_{-1},Y(v,z,\bz)]=\partial
Y(v,z,\bz),\quad [L_0,Y(v,z,\bz)]=z\partial Y(v,z,\bz)+Y(L_0
v,z,\bz),$ for any $v\in V,$ and similar conditions obtained by
replacing $L_n\ra \bL_n,$ $\partial\ra\bpartial.$ \item[(iv)]
$Y(v)\vac = e^{zL_{-1}+\bz\bL_{-1}}v$ for any $v\in V.$ \item[(v)]
$Y(v,z,\bz)Y(v',z',\bz')$ has only power-like singularities on the
diagonals $z=z'$ and $\bz=\bz'.$ \item[(vi)]
$Y(v,z,\bz)Y(v',z',\bz')-Y(v',z',\bz')Y(v,z,\bz)$ is a formal
distribution supported on the diagonal.
\end{itemize}

Recall that the Virasoro algebra is an infinite-dimensional Lie algebra spanned by elements $L_m,m\in\ZZ$ and the following
commutation relations:
$$
[L_m,L_n]=(m-n)L_{m+n}+c\frac{m^3-m}{12}\delta_{m,-n}.
$$
It is a unique central extension of the Witt algebra (the Lie algebra of vector fields on a circle).
The constant $c$ is called the central charge.
The Virasoro algebras spanned by $L_n$ and $\bL_n$ are called right-moving and left-moving, respectively.

There are certain variations of this definition. The modification
which we will need most amounts to replacing all spaces and maps
by their $\ZZ/2$-graded versions, and the ``commutativity'' axiom
(vi) with supercommutativity. From the physical viewpoint, this
means that we allow both fermions and bosons in our theory.
Another important property which must hold in any acceptable CFT
is the existence of a non-degenerate bilinear form on $V$ which is
compatible, in a suitable sense, with the rest of the data.
Finally, most CFTs of interest are ``left-right symmetric.'' This
means that exchanging $z$ and $\bz,$ and $L_n$ and $\bL_n,$ gives
an isomorphic CFT. We will only consider left-right symmetric
CFTs.

A more geometric approach to 2d CFT has been proposed by G.~Segal~\cite{Segal}. In Segal's approach, one starts with
a certain category
whose objects are finite ordered sets of circles, and morphisms are Riemann surfaces with oriented and analytically
parametrized boundaries.
Composition of morphisms is defined by sewing Riemann surfaces along boundaries with compatible orientations.
A 2d CFT is a projective functor from this category to the category of Hilbert spaces which satisfies certain
properties which are listed in ~\cite{Segal, Gaw}. (A projective functor from a category ${\mathcal C}$ to the category
of Hilbert spaces is the same as a functor from ${\mathcal C}$ to a category whose objects are Hilbert spaces, and
morphisms are equivalence classes of Hilbert space morphisms under the operation of multiplication by non-zero scalars.)
One can show that any 2d CFT in the sense of Segal's
definition gives rise to a 2d CFT in the sense of our algebraic definition (see e.g.~\cite{Gaw}). For example,
the vector space $V$ which appears in our algebraic definition is the Hilbert space associated to
a single circle in Segal's approach. The map $Y$ comes from considering the
morphism which corresponds to a Riemann sphere with three holes. Conversely, it appears that any ``algebraic''
2d CFT which is left-right symmetric and is equipped with a compatible inner product gives rise to a ``geometric''
2d CFT in genus zero (i.e. with Riemann surfaces restricted to have genus
zero).

An $N=1$ super-Virasoro algebra is a certain infinite-dimensional
Lie super-algebra which contains the ordinary Virasoro as a
subalgebra. Apart from the Virasoro generators $L_n,n\in \ZZ,$ it
contains odd generators $Q_n,n\in \ZZ.$ The additional commutation
relations read
$$
[L_m,Q_n]=\left(\frac{m}{2}-n\right) Q_{m+n},\quad [Q_m,Q_n]= \frac{1}{2} L_{m+n} + \frac{c}{12}m^2
\delta_{m,-n}.
$$
An $N=1$ superconformal field theory (SCFT) is a 2d CFT with an action on $V$ of two copies of
the $N=1$ super-Virasoro algebra which is compatible with other structures of the SCFT in a fairly obvious
sense.
\footnote{Strictly speaking, we are talking about the Ramond-Ramond sector of the SCFT here.}

An $N=2$ super-Virasoro algebra is a further generalization of the
Virasoro algebra. It is a certain infinite-dimensional Lie
super-algebra which contains the $N=1$ super-Virasoro as a
subalgebra. The even generators are $L_n,J_n,$ $n\in \ZZ.$ The odd
generators are $Q^\pm_n,n\in\ZZ.$ The commutators read,
schematically:
\begin{align}
&[L,L]\sim L,\quad  &[J,J]\sim central,\quad &[L,J]\sim J,\quad &[Q^\pm,Q^\pm]=0,\\
&[L,Q^\pm]\sim Q^\pm,\quad &[J,Q^\pm]\sim \pm Q^\pm,\quad &[Q^\pm,Q^\mp]\sim J+L+central.&
\end{align}
The precise form of the commutation relations ca be found in
~\cite{KO}. The $N=1$ super-Virasoro subalgebra is spanned by
$L_n$ and $Q_n=Q^+_n+Q^-_n.$ The relation between $N=1$ and $N=2$
super-Virasoro is analogous to the relation between the de Rham
and Dolbeault differentials on a complex manifold: $Q_n$ are
analogous to $d,$ while $Q^+_n$ and $Q^-_n$ are analogous to
$\partial$ and $\bpartial.$

An $N=2$ SCFT is an $N=1$ SCFT with an action of two copies of $N=2$
super-Virasoro algebra which is compatible with the remaining structures of the SCFT.
Thus we have a hierarchy of algebraic structures:
$$
{\rm Set\ of\ all\ CFTs} \supset {\rm Set\ of\ all\ N=1\ SCFTs} \supset {\rm Set\ of\ all\ N=2\ SCFTs}
$$
In fact, there is an even more general notion: 2d quantum field theory, without the adjective ``conformal.''
We will not discuss it in these lectures.

It is possible to give a definition of $N=1$ and $N=2$ superconformal field theories \`{a} la Segal.
The role of Riemann surfaces is played by 2d supermanifolds equipped with $N=1$ or $N=2$
superconformal structure.

An isomorphism of (super-)conformal field theories is a 1-1 map
$V\stackrel{\sim}{\to} V'$ which preserves all the relevant
structures. Two $N=2$ superconformal field theories can be
isomorphic as $N=1$ superconformal field theories without being
isomorphic as $N=2$ superconformal field theories. (When
physicists say that two (super-)conformal field theories are ``the
same'', they often neglect to specify which structures are
preserved by the isomorphism; this is usually clear from the
context.)

$N=2$ super-Virasoro algebra has an interesting automorphism called the mirror automorphism:
$$
\m: L_n\mapsto L_n,\quad J_n\mapsto -J_n,\quad Q^\pm_n\mapsto Q^\mp_n.
$$

Suppose we have a pair of $N=2$ superconformal field theories
which are isomorphic as $N=1$ SCFTs. Let $f:V\stackrel{\sim}{\to}
V'$ be an isomorphism. We say that $f$ is a (right) mirror
morphism of $N=2$ CFTs if it acts as the identity on the
``left-moving'' $N=2$ super-Virasoro, and acts by the mirror
automorphism on the ``right-moving'' $N=2$ super-Virasoro. This
makes sense because the mirror automorphism acts as the identity
on the $N=1$ super-Virasoro subalgebra of $N=2$ super-Virasoro
algebra. Exchanging left and right, we get the notion of a left
mirror morphism of $N=2$ CFTs. Finally, if $f$ acts as the mirror
automorphisms on both left and right super-Virasoro algebras, we
will say that $f$ is a target-space complex conjugation.

By definition, two $N=2$ SCFTs are mirror to each other if there
is a (left or right) mirror morphism between them. Clearly, if two
$N=2$ SCFTs are both (left-) mirror to a third $N=2$ SCFT, then
the first two SCFTs are isomorphic (as $N=2$ SCFTs). Thus the
mirror relation (for example, left) is an involutive relation on
the set of isomorphism classes of $N=2$ SCFTs. We stress that if
two $N=2$ SCFTs are mirror to each other, then they are isomorphic
as $N=1$ SCFTs, but usually not as $N=2$ SCFTs. Many explicit
examples of mirror pairs of $N=2$ SCFTs have been constructed in
the physics literature; one can construct more complicated
examples by means of tensor product, orbifolding, etc.

Now let us turn to the relation between $N=2$ SCFTs and Calabi-Yau manifolds. A physicist's Calabi-Yau
is a compact complex manifold with a trivial canonical class equipped with a K\"ahler class and
a B-field (an element of $H^2(X,\RR)/H^2(X,\ZZ)$).
It is believed that to any physicist's Calabi-Yau one can attach, in a natural way, an $N=2$ SCFT  which
depends on these geometric data.
One can give the following heuristic argument supporting the claim.
First of all, to any physicist's Calabi-Yau one can naturally attach
a classical field theory called the $N=2$ sigma-model. Its Lagrangian is given by an explicit, although rather
complicated, formula (see e.g.~\cite{Polchtext}). Infinitesimal symmetries of this classical field theory include
two copies of $N=2$
super-Virasoro algebra (with zero central charge). Second,
one can try to quantize this classical field theory while preserving $N=2$ superconformal invariance
(up to an unavoidable
central extension). The result of the quantization should be an $N=2$ SCFT.

Except for a few special cases, it is not known how to quantize
the sigma-model exactly. On the other hand, one has a perturbative
quantization procedure which works when the volume of the
Calabi-Yau is large. That is, if one rescales the metric by a
parameter $t\gg 1,$ $g_{\mu\nu}\ra t^2 g_{\mu\nu},$ and considers
the limit $t\ra \infty$ (so called large volume limit), then one
can quantize the sigma-model order by order in $1/t$ expansion. It
is believed that the resulting power series in $1/t$ has a
non-zero radius of convergence, and defines an actual $N=2$ SCFT.

It is natural to ask if it is possible to reconstruct a Calabi-Yau
starting from an $N=2$ SCFT; as we will discuss shortly, the
reconstruction problem does not have a unique answer. However,
some numerical characteristics of the ``parent'' Calabi-Yau $X$
can be determined rather easily. For example, the complex
dimension of $X$ is given by $c/3,$ where $c$ is the central
charge of the $N=2$ super-Virasoro algebra. One can also determine
the Hodge numbers $h^{p,q}(X)$ in the following manner.
Commutation relations of the $N=2$ super-Virasoro imply that the
operator $D_B=Q_0^++\bQ_0^+$ on $V$ squares to zero. $D_B$ is
known as a BRST operator (of type B, see below), and its
cohomology is called the BRST cohomology. The BRST cohomology
$\Ker\ D_B/\Im\ D_B$ is finite-dimensional in any reasonable $N=2$
SCFT. It is graded by the eigenvalues of the operators $J_0$ and
$\bJ_0$ known as left and right-moving R-charges. The Hodge number
$h^{p,q}(X)$ is simply the dimension of the component with
R-charges $p-\frac{n}{2}$ and $q-\frac{n}{2},$ where $n=\dim_\CC
X.$ This means, incidentally, that not every $N=2$ SCFT arises
from a Calabi-Yau manifold: those which do, must have integral
$n=c/3$ and integral spectrum of $J_0+\frac{n}{2}$ and
$\bJ_0+\frac{n}{2}$ in BRST cohomology. It is believed that any
$N=2$ SCFT with integral $c/3$ and integral spectrum of
$J_0+\frac{n}{2}$ and $\bJ_0+\frac{n}{2}$ is related to some
Calabi-Yau manifold, if one allows certain kinds of singular
Calabi-Yau manifolds, such as orbifolds.

One can show that if two Calabi-Yau are complex-conjugate, then the corresponding $N=2$ SCFTs are related by
target-space conjugation (in the sense explained above). This explains why the name ``target-space complex
conjugation'' was attached to a particular kind of morphisms of $N=2$ SCFTs. On the other hand, if
two Calabi-Yau manifolds
have $N=2$ SCFTs related by target-space complex conjugation, this does not imply that the Calabi-Yau manifolds
themselves are complex-conjugate; it merely implies that their $N=2$ SCFTs are related in a simple way.
Thus one obvious question is

Question 1. When do two Calabi-Yau manifolds produce isomorphic
$N=2$ SCFTs?

An answer to this question would interpret ``quantum symmetries'' of 2d SCFTs in geometric terms. Another question of this
kind is

Question 2. When do two Calabi-Yau manifolds produce $N=2$ SCFTs
which are mirror to each other?

We say that two Calabi-Yau manifolds are related by mirror symmetry if the corresponding $N=2$ SCFTs are mirror.

First non-trivial examples of mirror pairs of Calabi-Yau manifolds
have been constructed by B. Greene and R. Plesser~\cite{GP}. The
simplest example in complex dimension three (this dimension is the
most interesting one from the physical viewpoint) is the
following: one of the Calabi-Yau manifolds is the Fermat quintic
$x^5+y^5+z^5+v^5+w^5=0$ in $\CC\PP^4,$ while the other one is
obtained by taking a quotient of the Fermat quintic by a certain
action of $(\ZZ/5)^3$ and blowing up the fixed points. We will not
try to explain why these two Calabi-Yau manifolds are mirror. (The
original argument~\cite{GP} relied on a conjectural equivalence
between the $N=2$ SCFT corresponding to the Fermat quintic and a
certain integrable $N=2$ SCFT constructed by D. Gepner. Later this
issue has been studied in detail by E. Witten~\cite{phases} and
now has the status of a physical ``theorem.'')

The answer to the first question  is highly non-trivial. This can
be seen already in the case when $X$ is a complex torus with a
flat metric (Lecture 2). For example, the torus and its dual give
the same $N=2$ SCFT, even though they are usually not isomorphic
as complex manifolds. The answer to Question 2 -- characterization
of the mirror relation in geometric terms -- is the ultimate goal
of the Mirror Symmetry program.

On the most basic level, a mirror relation between $X$ and $X'$
implies a relation between their Hodge numbers $h^{p,q}(X)$ and
$h^{p,q}(X').$ To see how this comes about, note that along with
the cohomology of $D_B=Q^+_0+\bQ^+_0$ we may also consider the
cohomology of $D_A=Q^-_0+\bQ^+_0.$ It can be shown that in any
$N=2$ SCFT satisfying the integrality constraint these two
cohomologies are isomorphic as bi-graded vector spaces~\cite{LVW}.
Now note that if $X$ is mirror to $X',$ then the cohomology of
$D_A(X)$ in bi-degree $(\alpha,\beta)$ is isomorphic to cohomology
of $D_B(X')$ in bi-degree $(-\alpha,\beta).$ Recalling the
relation between the Hodge numbers of $X$ and cohomology of
$D_B(X),$ we infer an important relation
\begin{equation}\label{Hodge}
h^{p,q}(X)=h^{n-p,q}(X').
\end{equation}
If one plots the Hodge numbers of a complex manifold on a plane
with coordinates $p-q$ and $p+q-n,$ the resulting table has the
shape of a diamond (the Hodge diamond). For any Calabi-Yau
manifold the Hodge diamond is unchanged by a rotation by
$180^\circ$ degrees. The relation~(\ref{Hodge}) means that the
Hodge diamonds of mirror Calabi-Yau manifolds are related by a
rotation by $90^\circ$ degrees.

Of course, the existence of a mirror relation between $X$ and $X'$ implies much more than this. The most promising
approach to the problem
of characterizing the mirror relation in geometric terms has been proposed in 1994 by M.~Kontsevich.
In the remainder of this lecture
we will sketch Kontsevich's proposal and its interpretation in physical terms.

A physicist's Calabi-Yau has both a complex structure and a
symplectic structure (the K\"ahler form). One can gain a
considerable insight into the Mirror Symmetry Phenomenon by
focusing one of the two structures. More precisely, if the B-field
is present, we combine the K\"ahler form $\omega$ and the $(1,1)$
part of the B-field into a ``complexified K\"ahler form.'' We will
regard the latter as parametrizing an ``extended symplectic moduli
space'' of $X.$ Similarly, we regard the $(0,2)$ part of the
B-field and the complex structure moduli as parametrizing an
``extended complex structure moduli space'' of $X.$ We would like
to isolate some aspects of the $N=2$ SCFT which depend either only
on the extended complex moduli, or only on the extended symplectic
moduli. The procedure for doing this was proposed by E.
Witten~\cite{topsigma,Wittentwist} and is known as topological
twisting.

Witten's construction rests on the observation that many $N=2$ SCFTs have finite-dimensional sectors
which are topological field theories, i.e. do not depend on the 2d metric (the metric on the world-sheet, if we use
string theory terminology.) In fact, for many $N=2$ SCFTs there are two such sub-theories ; they are known as
A- and B-models. $N=2$ SCFTs related to Calabi-Yau manifolds belong precisely to this class of
theories.

Let us recall some basic facts about 2d topological field theories. These theories are similar to, but much simpler than,
2d CFTs. They can be described by axioms similar to Segal's axioms~\cite{Atiyah}.
One starts with a category whose objects
are finite ordered sets of oriented and parametrized circles and morphisms are oriented 2d manifolds (without complex
structure) bounding the circles.
A 2d topological field theory is a functor from this category to
the category of finite-dimensional (graded) vector spaces which satisfies certain requirements similar to Segal's.
As for 2d CFTs, there is a purely algebraic reformulation of this definition. It turns out that the ``topological'' counterpart
of the notion of a conformal field theory is the well-known notion of a super-commutative Frobenius algebra, i.e. a
super-commutative algebra with an invariant metric (see e.g.~\cite{Dijkgraaf}).

A detailed discussion of Witten's procedure for constructing a 2d
TFT out of an $N=2$ SCFT is beyond the scope of these lectures.
Roughly speaking, one passes from the space $V$ to its BRST
cohomology. One can show that the state-operator correspondence
$Y$ descends to a super-commutative algebra structure on the BRST
cohomology.  The invariant metric on BRST cohomology comes from a
metric on $V.$

Note that we have two essentially different choices of a BRST
operator: $D_A$ or $D_B.$ (One can also consider
$D_A'=Q^+_0+\bQ^-_0$ and $D_B'=Q^-_0+\bQ^-_0,$ but these can be
trivially related to $D_A$ and $D_B$ by replacing $X$ with its
complex-conjugate.) Thus Witten's construction associates to any
physicist's Calabi-Yau $X$ two topological field theories, called
the A-model and the B-model, respectively.

It turns out that the A-model does not change as one varies the
extended complex structure moduli, while the B-model does not
depend on the extended symplectic moduli. In other words, the
A-model isolates the symplectic aspects of the Calabi-Yau, while
the B-models isolates the complex ones. In fact, the state space
of the A-model is naturally isomorphic to the de Rham cohomology
$H^{*,*}(X),$ while the state space of the B-model is naturally
isomorphic to the Dolbeault cohomology
$$
H^*(\Lambda^* T^{1,0}X),
$$
where $T^{1,0}X$ is the holomorphic tangent bundle of $X.$ For a
Calabi-Yau manifold, $H^q(\Lambda^p T^{1,0}X)$ is isomorphic to
$H^{n-p,q}(X),$ but not canonically: the isomorphism depends on
the choice of a holomorphic section of the canonical line bundle.
Note also that the spaces of the A and B-models are bi-graded.
From the physical point of view, the bi-grading comes from the
decomposition of the state spaces into the eigenspaces of $J_0$
and $\bJ_0.$

The algebra structure in the B-case is the obvious one, while in
the A-case it is a deformation of the obvious one which depends on
the extended symplectic structure on $X$; $H^*(X)$ equipped with
this deformed algebra structure is known as the quantum cohomology
ring of $X.$

Mirror symmetry acts on A and B-models in a very simple way. It is
easy to see that the mirror automorphism exchanges $D_A$ and
$D_B.$ Thus if $X$ and $X'$ are a mirror pair of physicist's
Calabi-Yau manifolds, then the A-model of $X$ is isomorphic to the
B-model of $X',$ and vice versa. We will say that $X$ and $X'$ are
weakly mirror if the A-model of $X$ is isomorphic to the B-model
of $X',$ and vice-versa. The notion of weak mirror symmetry is
easier to work with, since one can define A and B-models of a
Calabi-Yau directly, without appealing to the ill-defined
quantization of the sigma-model. But clearly a lot of information
is lost in the course of topological twisting, and one would like
to find some richer objects associated to an $N=2$ SCFT.

M.~Kontsevich proposed that a suitable enriched version of the
B-model is the bounded derived category of coherent sheaves on
$X,$ which we will denote $\bD^b(X),$ while the enriched version
of the A-model is some version of the derived Fukaya category of
$X,$ which will be denoted $\bD\cF_0(X).$ In other words, he
conjectured that if $X$ and $X'$ are mirror, then $\bD^b(X)$ is
equivalent to $\bD\cF_0(X')$ and vice versa. This is known as the
Homological Mirror Conjecture (HMC). If the converse statement
were true, this would completely answer Question 1 and 2.

Let us sketch the definitions of these two categories. Let $X$ be
a smooth complex manifold. An object of the bounded derived
category on $X$ is a bounded complex of holomorphic vector bundles
on $X,$ i.e. a finite sequence of holomorphic vector bundles and
morphisms between them
$$
0\ra\cdots\ra E_{n-1} \ra E_n \ra E_{n+1}\ra \cdots \ra 0,
$$
so that the composition of any two successive morphisms is zero.
We remind the reader that the cohomology of such a complex is a
sequence of coherent sheaves on $X.$ To define morphisms in the
derived category, we first consider the category of bounded
complexes ${\bf C}^b(X),$ where morphisms are defined as chain
maps between complexes. A morphism in this category is called a
quasi-isomorphism if it induces an isomorphism on the cohomology
of complexes. The idea of the derived category is to identify all
quasi-isomorphic complexes. That is, the bounded derived category
$\bD^b(X)$ is obtained from ${\bf C}^b(X)$ by formally inverting
all quasi-isomorphisms. In this definition, one can replace
holomorphic vector bundles by arbitrary coherent sheaves; the
resulting derived category is unchanged. Lecture 3 will discuss
derived categories in more detail.

While coherent sheaves and their complexes are very familiar
creatures and are the basic tool of algebraic geometry, the
derived Fukaya category $\bD\cF_0(X)$ is a much more recent
invention. It is obtained by a rather complicated algebraic
procedure from a certain geometrically defined category called the
Fukaya category $\cF(X).$ The latter has been introduced by
K.~Fukaya in \cite{Fukaya}. Actually, $\cF(X)$ is not quite a
category: there are additional structures on morphisms (multiple
compositions and the differential), and the composition of
morphisms is associative only up to a chain homotopy. Such a
structure is called an $A_\infty$ category. The Fukaya category
depends only on the extended symplectic structure on $X.$ Objects
of the Fukaya category are, roughly speaking, triples
$(Y,E,\nabla),$ where $Y$ is a Lagrangian submanifold of $X,$ $E$
is a complex vector bundle on $Y$ with a Hermitian metric, and
$\nabla$ is a flat unitary connection on $E.$ This definition is
only approximate, for several reasons. First of all, not every
Lagrangian submanifold $Y$ is allowed: the so-called Maslov class
of $Y$ must vanish (the Maslov class is a class in $H^1(Y,\ZZ)$).
Second, $Y$ has to be a graded Lagrangian submanifold (this notion
was defined in 1968 by J.~Leray for the case of Lagrangian
submanifolds in a symplectic vector space, and generalized by
Kontsevich to the Calabi-Yau case). Third, it is not completely
clear if the flat connection $\nabla$ has to be unitary. There are
some indications that it might be necessary to relax this
condition and require instead the eigenvalues of the monodromy
representation to have unit absolute value. Fourth, in the
presence of the B-field $\nabla$ must be projectively flat rather
than flat~\cite{KO}.

The space of morphisms in the Fukaya category is defined by means of the Floer complex. This will be
discussed in Lecture 2.

The relation between the Homological Mirror Conjecture and the A
and B-models is the following~\cite{Konts}. Given a triangulated
category (more precisely, an ``enhanced triangulated
category''~\cite{BK}), one can study its deformations. Information
about deformations is encoded in the Hochschild cohomology of the
category in question. In the case of the derived category of
coherent sheaves, Hochschild cohomology seems to coincide with the
 cohomology of the exterior algebra of the holomorphic tangent
bundle, i.e. the state space of the B-model.\footnote{We say
``seems'', because there is no complete proof of this yet.}
Kontsevich conjectured that the Hochschild cohomology of the
derived Fukaya category is the quantum cohomology ring of $X,$
i.e. the state space of the A-model~\cite{Konts}. Thus the
equivalence of $\bD^b(X)$ and $\bD\cF_0(X')$ is likely to imply
that the B-model of $X$ is isomorphic (as a 2d TFT) to the A-model
of $X'.$ In other words, Homological Mirror Symmetry seems to
imply weak mirror symmetry.

Homological Mirror Symmetry Conjecture also has a clear physical
meaning. The physical idea is to generalize the notion of a 2d TFT
to allow the 2d world-sheet to have boundaries (see
e.g.~\cite{WittenCS}). This generalization also makes sense in the
full $N=2$ SCFT and leads to the notion of a D-brane, which plays
a very important role in string theory~\cite{Polch}. A D-brane is
a nice boundary condition for the SCFT. It is not completely clear
what this means in the quantum case, so let us discuss this notion
using classical field theory. A classical 2d field theory is
defined by an action which is an integral of a local Lagrangian
over the 2d world-sheet. So far we took the world-sheet to be a
cylinder whose noncompact direction was parametrized by the
``time'' variable. Thus the space was topologically a circle. Now
let us take the space to be an interval $I.$ The world-sheet
becomes $\RR\times I.$ In order for the classical field theory to
be well-defined, we require the the Cauchy problem for the
Euler-Lagrange equations to have a unique solution, at least
locally. This requires imposing a suitable boundary condition on
the fields and their derivatives on the boundary of the
world-sheet. For example, one can impose Dirichlet boundary
conditions (i.e. vanishing) on some scalar fields which appear in
the Lagrangian.\footnote{The letter D in the word ``D-brane''
actually came from ``Dirichlet.''} A classical D-brane is simply a
choice of such a boundary condition.

If the classical field theory has some symmetries, it
is reasonable to require the boundary condition to preserve this symmetry. For example, $N=1$ sigma-models
have two copies of $N=1$ super-Virasoro algebra as their symmetries. It is not possible to preserve both of them,
but there exist many boundary conditions which preserve the diagonal subalgebra. Such boundary conditions are ordinary
D-branes of superstring theory~\cite{Polch}. In the $N=2$ case, we have two copies of $N=2$ super-Virasoro, and we may require
the boundary condition to preserve the diagonal $N=2$ super-Virasoro. Such boundary conditions are called D-branes
of type B, or simply B-branes, because they are related to the B-model (see below). One can also exploit the existence of the
mirror automorphism $\m$ and consider boundary conditions which preserve a different $N=2$ super-Virasoro subalgebra,
namely the one spanned by
$$
L_n+\bL_n,\quad -J_n+\bJ_n,\quad Q^-_n+\bQ^+_n, \quad Q^+_n+\bQ^-_n,\quad n\in\ZZ.
$$
The corresponding branes are called D-branes of type A, or simply A-branes.

Given a classical D-brane, we can try to quantize a classical
field theory on $\RR\times I$ with boundaries, which related to
this $D$-brane, so that the quantized theory has one copy of $N=2$
super-Virasoro as its symmetry algebra. If such a quantization is
possible, we say that the classical D-brane is quantizable, and
the classical D-brane together with its quantization will be
called a quantum D-brane. This is not a very satisfactory way to
define quantum D-branes, and it remains an interesting problem to
find a satisfactory and fully quantum definition of a boundary
condition for a 2d SCFT.

Now let us turn to the relation of A and B-branes with A and
B-models. A and B-models are obtained from the $N=2$ SCFT by
topological twisting. Roughly speaking, twisting amounts to
truncating the theory to the cohomology of $D_A$ or $D_B.$ Now
note that $D_A$ (resp. $D_B$) sits in the $N=2$ super-Virasoro
which is preserved by the A-type (resp. B-type) boundary
condition. The significance of this is that the A-twist is
consistent with A-type boundary conditions, while the B-twist is
consistent with B-type boundary conditions. Thus an A-brane (resp.
B-brane) gives rise to a consistent boundary condition for the
A-model (resp. B-model).\footnote{Axioms for boundary conditions
in 2d TFTs have been discussed by G.~Moore and G.~Segal~\cite{MS}
and C.I.~Lazaroiu~\cite{Laz}.}

One can show that the set of A-branes (or B-branes) has the
structure of a category. The space of morphisms between two branes
$A$ and $A'$ is simply the space of states of the 2d TFT on
$\RR\times I,$ with boundary conditions on the two ends
corresponding to $A$ and $A'.$ Equivalently, one considers the
state space of the $N=2$ SCFT on an interval, and computes its
BRST cohomology with respect to $D_A.$ The composition of
morphisms can be defined with the help of the state-operator
correspondence $Y$ (or rather, its analogue in the case of a 2d
SCFT with boundaries).

To summarize, to any physicist's Calabi-Yau we can attach two
categories: the categories of A-branes and B-branes. One can argue
that the category of A-branes (resp. B-branes) does not depend on
the extended complex (resp. extended symplectic)
moduli~\cite{WittenCS}. One can think of these categories as the
enriched versions of the A and B-models: while the A-model is a
TFT on a world-sheet without boundaries, the totality of A-branes
corresponds to the same 2d TFT on a world-sheet with boundaries
and with all possible boundary conditions. Further, it is obvious
that if two Calabi-Yau manifolds are related by a mirror morphism,
then the A-brane category of the first manifold is equivalent to
the B-brane category of the second one, and vice versa. Obviously,
if two $N=2$ SCFTs related two Calabi-Yau manifolds are
isomorphic, then the corresponding categories of A-branes (and
B-branes) are simply equivalent.

The Homological Mirror Conjecture would follow if we could prove
that the category of A-branes (resp. B-branes) is equivalent to
$\bD\cF_0(X)$ (resp. $\bD^b(X)$). Alas, we cannot hope to prove
this, because we do not have an honest definition of a D-brane.
What we do know is that holomorphic vector bundles are examples of
B-branes, and objects of the Fukaya category are examples of
A-branes~\cite{WittenCS}. Furthermore, Witten showed that in this
special case morphisms in the category of B-branes and A-branes
agree with morphisms in the categories $\bD^b(X)$ and
$\bD\cF_0(X),$ respectively~\cite{WittenCS}. This computation
served as a motivation for Kontsevich's conjecture. More recently
it was shown that more general coherent sheaves, as well as
complexes of coherent sheaves, are also valid B-branes. On the
other hand, it has been shown recently that the Fukaya category is
only a full subcategory of the category of A-branes, that is, for
some $X$ there exist A-branes which are not isomorphic to any
object of $\bD\cF_0(X)$~\cite{KO2}. This means that the symplectic
side of Kontsevich's conjecture needs substantial modification.
This will be discussed in more detail in Lecture 4.

From the mathematical point of view, the Homological Mirror
Conjecture is not well-defined, because it is not clear how to
quantize the sigma-model for an arbitrary Calabi-Yau manifold. But
there is a class of Calabi-Yaus for which the sigma-model can be
quantized, and the corresponding $N=2$ SCFTs can be described
quite explicitly. These are complex tori with a flat K\"ahler
metric and arbitrary B-field. One can easily determine which pairs
of such Calabi-Yaus give mirror $N=2$ SCFTs; the resulting
criterion can be expressed in terms of linear algebra~\cite{KO}.
(This will be discussed in Lecture 2). Thus in the case of flat
tori the Homological Mirror Conjecture is mathematically
well-defined (modulo the issues related to the precise definition
of the Fukaya category). In ~\cite{PZ} A. Polishchuk and E. Zaslow
proved it for tori of real dimension two. But there are strong
arguments showing that for higher-dimensional tori the Homological
Mirror Conjecture cannot hold, unless one substantially enlarges
the Fukaya category by adding new objects. This will be discussed
in Lecture 4.

For more general Calabi-Yaus, one can take the Homological Mirror Conjecture
as an attempt to give a {\it mathematical definition} of the mirror relation.
Then the main issues are the precise definition of the Fukaya category,
and the verification that the numerous mirror pairs proposed by physicists and mathematicians are
in fact mirror  in the sense of the Homological Mirror Conjecture.

Another approach to mirror symmetry has been proposed
in~\cite{SYZ} and is known as the SYZ Conjecture. According to
this conjecture, mirror Calabi-Yau manifolds admit fibrations by
special Lagrangian tori which in some sense are ``dual'' to each
other. Recently a relation between the Homological Mirror
Conjecture and the SYZ Conjecture has been studied in the paper
~\cite{KS}.

\section{Mirror Symmetry for Flat Complex Tori}

In this lecture we describe $N=2$ superconformal field theories
(SCFT) related to complex tori $T$ endowed with a flat K\"ahler
metric $G$ and a constant 2-form $B$ (the B-field). We will give a
criterion when two such data $(T, G, B)$ and $(T', G', B')$
produce isomorphic $N=2$ SCFT and when they produce $N=2$ SCFT
which are mirror to each other.

As explained in Lecture 1, to define an $N=2$ superconformal field
theory we need to specify an infinite-dimensional $\ZZ_2$-graded
vector space of states $V,$ a vacuum vector $|vac\rangle,$ a
state-operator correspondence  $Y$ from $V$ to the space of
``formal fraction power series in $z, \bar{z}$ with coefficients
in $\End(V)$'' and, finally,  the super-Virasoro elements $L,
\bar{L}, Q^{\pm}, \bar{Q}^{\pm}, J, \bar{J}$ which enter into the
definition of the superconformal structure (see Lecture 1).

We start with some notations. Let  $\Gamma\cong \ZZ^{2d}$ be a
lattice in a real vector space $U$ of dimension $2d,$ and let
$\Gamma^*\subset U^*$ be the dual lattice. Consider real tori
$T=U/\Gamma,\ T^*=U^*/\Gamma^*.$ Let $G$ be a metric on $U,$ i.e.
a positive symmetric bilinear form on $U,$ and let $\B$ be a real
skew-symmetric bilinear form on $U.$ Denote by $l$  the natural
pairing $\Gamma\times \Gamma^*\to\ZZ.$ (The natural pairing
$U\times U^*\to\RR$ will be also denoted as $l.$) Choose
generators $e_1,\ldots,e_{2d}\in \Gamma.$ The components of an
element $w\in \Gamma$ in this basis will be denoted by
$w^i,i=1,\ldots,2d.$ The components of an element $m\in\Gamma^*$
in the dual basis will be denoted by $m_i, i=1,\ldots,2d.$ We also
denote by $G_{ij},\ \B_{ij}$ the components of $G,\ \B$ in this
basis. It will be apparent that the superconformal field theory
which we construct does not depend on the choice of basis in
$\Gamma.$ In the physics literature $\Gamma$ is sometimes referred
to as the winding lattice, while $\Gamma^*$ is called the momentum
lattice.

Consider a triple $(T,G,\B).$ With any such triple we associate an
$N=2$ superconformal field theory $\cV$ which may be regarded as a
quantization of the supersymmetric sigma-model.

The state space of the SCFT $\cV$ is
$$
V=\cH_b\ot_\CC\cH_f\ot_\CC \grp.
$$
Here $\cH_b$ and $\cH_f$ are bosonic and fermionic Fock spaces
defined below, while $\grp$ is the space of the group algebra of
$\Gamma\op\Gamma^*$ over $\CC.$

To define $\cH_b,$ consider an algebra over $\CC$ with generators
$\al^i_s,\bal^i_s,\ i=1,\ldots,2d, s\in \ZZ\backslash 0$ and
relations
\begin{equation}\label{CCR}
[\al^i_s,\al^j_p]=s\Gi^{ij}\delta_{s,-p}, \quad
[\bal^i_s,\bal^j_p]=s\Gi^{ij}\delta_{s,-p}, \quad
[\al^i_s,\bal^j_p]=0.
\end{equation}
If $s$ is a positive integer, $\al^i_{-s}$ and $\bal^i_{-s}$ are
called left and right bosonic creators, respectively, otherwise
they are called left and right bosonic annihilators. Either
creators or annihilators are referred to as oscillators.

The space $\cH_b$ is defined as the space of polynomials of even
variables $a^i_{-s},\ba^i_{-s}, i=1,\ldots,2d, s=1,2,\ldots.$ The
bosonic oscillator algebra acts on the space $\cH_b$ via
$$
\begin{array}{lllllll}
\al^i_{-s}&\mapsto& a^i_{-s}\cdot, &\qquad& \bal^i_{-s}&\mapsto&
\ba^i_{-s}\cdot, \\
\al^i_s &\mapsto&\ds s\Gi^{ij}\frac{\partial}{\partial a^j_{-s}},
& \qquad& \bal^i_s&\mapsto&\ds s\Gi^{ij}\frac{\partial}{\partial
\ba^j_{-s}},
\end{array}
$$
for all positive $s.$ This is the Fock-Bargmann representation of
the bosonic oscillator algebra. The vector $1\in \cH_b$ is
annihilated by all bosonic annihilators and will be denoted
$|vac_b\rangle.$

The space $\cH_b$ will be regarded as a $\ZZ_2$-graded vector
space with a trivial (purely even) grading. It is clear that
$\cH_b$ can be decomposed as $\fH_b\ot \bar{\fH}_b,$ where $\fH_b$
(resp. $\bar{\fH}_b$) is the bosonic Fock space defined using only
the left (right) bosonic oscillators.

To define $\cH_f,$ consider an algebra over $\CC$ with generators
$\psi^i_s,\tpsi^i_s,\ i=1,\ldots,2d, s\in\ZZ+\frac{1}{2}$ subject
to relations
\begin{equation}\label{CAR}
\{\psi^i_s,\psi^j_p\}=\Gi^{ij}\delta_{s,-p}, \quad
\{\bpsi^i_s,\bpsi^j_p\}=\Gi^{ij}\delta_{s,-p}, \quad
\{\psi^i_s,\bpsi^j_p\}=0.
\end{equation}
If $s$ is positive, $\psi^i_{-s}$ and $\tpsi^i_{-s}$ are called
left and right fermionic creators respectively, otherwise they are
called left and right fermionic annihilators. Collectively they
are referred to as fermionic oscillators.

The space $\cH_f$ is defined as the space of skew-polynomials of
odd variables $\theta^i_{-s},\btheta^i_{-s}, i=1,\ldots,2d,
s=1/2,3/2,\ldots.$ The  fermionic oscillator algebra~(\ref{CAR})
acts on $\cH_f$ via
$$
\begin{array}{lllllll}
\psi^i_{-s}&\mapsto& \theta^i_{-s}\cdot, &
\qquad& \bpsi^i_{-s}&\mapsto& \btheta^i_{-s}\cdot,\\
\psi^i_s&\mapsto&\ds \left(G^{-1}\right)^{ij}\frac{\partial}
{\partial \theta^j_{-s}}, & \qquad& \bpsi^i_s&\mapsto&
\ds\left(G^{-1}\right)^{ij}\frac{\partial}{\partial
\btheta^j_{-s}},
\end{array}
$$
for all positive $s\in\ZZ+\frac{1}{2}.$ This is the Fock-Bargmann
representation of the fermionic oscillator algebra. The vector
$1\in \cH_f$ is annihilated by all fermionic annihilators and will
be denoted $|vac_f\rangle.$ The fermionic Fock space has a natural
$\ZZ_2$ grading such that $|vac_f\rangle$ is even. It can be
decomposed as $\fH_f\ot \bar{\fH}_f,$ where $\fH_f$ (resp.
$\bar{\fH}_f$) is constructed using only the left (right)
fermionic oscillators.

For $w\in \Gamma,\ m\in \Gamma^*$ we will denote the vector $w\op
m\in \grp$ by $(w,m).$ We will also use a shorthand $|vac,
w,m\rangle$ for
$$
|vac_b\rangle\ot|vac_f\rangle\ot (w,m).
$$

To define $\cV,$ we have to specify the vacuum vector, $T,\bT,$
and the state-operator correspondence $Y.$ But first we need to
define some auxiliary objects. We define the operators $\hw:V\to
V\ot\Gamma$ and $\hm:V\to V\ot\Gamma^*$ as follows:
$$
\hw^i: b\ot f\ot (w,m)\mapsto w^i(b\ot f\ot(w,m)), \quad \hm_i:
b\ot f\ot(w,m)\mapsto m_i(b\ot f\ot (w,m)).
$$
We also set
$$
\partial X^j(z)  = \frac{1}{z} \Gi^{jk}\hp_k +
\sideset{}{'}\sum_{s=-\infty}^\infty \frac{\al_s^j}{
z^{s+1}},\qquad
 \bpartial X^j(\bz)  =  \frac{1}{\bz}
\Gi^{jk}\htp_k +\sideset{}{'}\sum_{s=-\infty}^\infty
\frac{\bal_s^j}{\bz^{s+1}},
$$
$$
\psi^j(z) = \sum_{r\in \ZZ+\frac{1}{2}}
\frac{\psi^j_r}{z^{r+1/2}}, \qquad \tpsi^j(\bz) = \sum_{r\in
\ZZ+\frac{1}{2}} \frac{\tpsi^j_r}{\bz^{r+1/2}},
$$
where a prime on a sum over $s$ means that the term with $s=0$ is
omitted, and $\hp_k$ and $\htp_k$ are defined by
$$
\hp_k=\frac{1}{\sqrt{2}}(\hm_k+\left(-\B_{kj}-G_{kj}\right)\hw^j),\quad
\htp_k=\frac{1}{\sqrt{2}}(\hm_k+\left(-\B_{kj}+G_{kj}\right)\hw^j).
$$
Note that we did not define $X^j(z,\bz)$ themselves, but only
their derivatives. The reason is that the would-be field
$X^j(z,\bz)$ contains terms proportional to $\log  z$ and $\log
\bz,$ and therefore is not a ``fractional power series.''

The vacuum vector of $\cV$ is defined by
$$
|vac\rangle=|vac,0,0\rangle.
$$

The general formula  for the state-operator correspondence $Y$ is
complicated and can be found in \cite{KO}. We will only list a few
special cases of the state-operator correspondence. The states
$\al^j_{-s}|vac,0,0\rangle$ and $\bal^j_{-s}|vac,0,0\rangle$ are
mapped by $Y$ to
$$
\frac{1}{(s-1)!}\partial^s X^j(z), \qquad
\frac{1}{(s-1)!}\bpartial^s X^j(\bz).
$$
The states $\psi^j_{-s}|vac,0,0\rangle$ and
$\tpsi^j_{-s}|vac,0,0\rangle$ are mapped to
$$
\frac{1}{(s-\frac{1}{2})!}\partial^{s-1/2}\psi^j(z), \qquad
\frac{1}{(s-\frac{1}{2})!}\bpartial^{s-1/2}\tpsi^j(\bz).
$$

\bigskip


To define an $N=2$ superconformal structure on $\cV,$ we need to
choose a complex structure $I$ on $U$ with respect to which $G$ is
a K\"ahler metric. Let $\omega=GI$ be the corresponding K\"ahler
form. Then the left-moving vectors are defined as follows:
\begin{align*}
L&=\frac{1}{2} G\left(a_{-1},a_{-1}\right)-
\frac{1}{2}G\left(\theta_{-1/2},\theta_{-3/2}\right),\\
Q^\pm&=\frac{-i}{4\sqrt 2}\left(G\mp
i\omega\right)\left(\theta_{-1/2},a_{-1}\right),\\
J&=-\frac{i}{2}\omega\left(\theta_{-1/2},\theta_{-1/2}\right).
\end{align*}
The right-moving vectors $\bar{L},$ $\bar{Q}^\pm$ and $\bar{J}$
are defined by the same expressions with $a$ replaced by $\bar{a}$
and $\theta$ replaced by $\bar{\theta}.$

{\it An isomorphism} of $N=2$ SCFTs  is an isomorphism of
underlying vector spaces $ f: V\stackrel{\sim}{\longrightarrow}
V', $ which preserves the state-operator correspondence $
Y'(f(a))f(b)=f(Y(a)b) $ and acts on the generators of both left
and right  super-Virasoro algebras as the identity map.

{\it A mirror morphism} between two $N=2$ SCFTs is an isomorphism
between the underlying $N=1$ SCFTs which induces the following map
on the generators of left/right super-Virasoro algebras:
\begin{align*}
&f(L)= L',\ f(Q^+)=Q^{-'},\ f(Q^-)=Q^{+'},\ f(J)=-J',\\
&f(\bL)=\bL',\ f(\bQ^+)=\bQ^{+'},\ f(\bQ^-)=\bQ^{-'},\
f(\bJ)=\bJ'.
\end{align*}
A composition of two mirror morphisms is an isomorphism of $N=2$
SCFTs.

Now we can describe when two different quadruples
$(\Gamma,I,G,\B)$ and $(\Gamma',I',G',\B')$ yield isomorphic $N=2$
SCFTs and when they are mirror symmetric.

The natural pairing $l:\Gamma\op\Gamma^*\ra \ZZ$ induces  a
natural $\ZZ$-valued symmetric bilinear form $q$ on
$\Gamma\op\Gamma^*$ defined by
\begin{equation}\label{bf}
q((w_1,m_1),(w_2,m_2))=l(w_1,m_2)+l(w_2,m_1), \qquad
w_{1,2}\in\Gamma,\ m_{1,2}\in\Gamma^*.
\end{equation}
Given $G,I,\B,$ we can define two complex structures on $T\times
T^*$:
\begin{align}
\cI(I,\B)&=\begin{pmatrix} I & 0 \\ \B I+I^t\B & -I^t \end{pmatrix}, \\
\cJ(G,I,\B)&=\begin{pmatrix} -IG^{-1}\B & IG^{-1} \\ GI-\B
IG^{-1}\B & \B IG^{-1}
\end{pmatrix}.
\end{align}
The notation here is as follows. We regard $\cI$ and $\cJ$ as
endomorphisms of $U\op U^*,$ and write the corresponding matrices
in the basis in which the first $2d$ vector span $U,$ while the
remaining vectors span $U^*.$ In addition, $G$ and $B$ are
regarded as elements of $\Hom(U,U^*),$ and $I^t$ denotes the
endomorphism of $U^*$ conjugate to $I.$

It is easy to see that $\cJ$ depends on $G,I$ only in the
combination $\omega=GI,$ i.e. it depends only on the symplectic
structure on $T$ and the B-field. There is also a third natural
complex structure $\tilde{\cI}$ on $T\times T^*,$ which is simply
the complex structure that $T\times T^*$ gets because it is a
Cartesian product of two complex manifolds:
$$
\tilde{\cI}=\begin{pmatrix} I & 0 \\ 0 & -I^t \end{pmatrix}.
$$
This complex structure will play only a minor role in what
follows. Note that $\cI$ coincides with $\tilde{\cI}$ if and only
if $\B^{(0,2)}=0.$

\begin{theorem}\label{ntwoiso}\cite{KO}
$\SCFT(\Gamma,I,G,\B)$ is isomorphic to $\SCFT(\Gamma',I',G',\B')$
if and only if there exists an isomorphism of lattices
$\Gamma\op\Gamma^*$ and $\Gamma'\op\Gamma^{'*}$ which takes $q$ to
$q',$ $\cI$ to $\cI',$ and $\cJ$ to $\cJ'.$
\end{theorem}

The second theorem describes when $(T,I,G,\B)$ is mirror symmetric
to $(T',I',G',\B').$
\begin{theorem}\label{ntwomirror}\cite{KO}
$\SCFT(\Gamma,I,G,\B)$ is mirror to $\SCFT(\Gamma',I',G',\B')$ if
and only if there is an isomorphism of lattices
$\Gamma\op\Gamma^*$ and $\Gamma'\op\Gamma^{'*}$ which takes $q$ to
$q',$ $\cI$ to $\cJ',$ and $\cJ$ to $\cI'.$
\end{theorem}

This theorem allows to give many examples of mirror pairs of tori.
Suppose that we are given a complex torus $(T,I)$ with a constant
K\"ahler form $\omega,$ and suppose that $T=A\times B,$ where $A$
and $B$ are Lagrangian sub-tori with respect to $\omega.$ In
particular, the lattice $\Gamma$ decomposes as
$\Gamma_A\op\Gamma_B.$ Let $\hA$ be the dual torus for $A,$ and
let $T'=\hA\times B.$ The lattice corresponding to $T'$ is
$\Gamma'=\Gamma_A^*\op \Gamma_B,$ and there is an obvious
isomorphism from $\Gamma\op \Gamma^*$ to $\Gamma'\op\Gamma^{'*}$
which takes $q$ to $q'.$ We let $\cI'$ and $\cJ'$ to be the image
of $\cJ$ and $\cI,$ respectively, and invert the relationship
between $\cI',\cJ'$ and $I',\omega',B'$ to find the complex
structure, the symplectic form, and the B-field on $T'.$ It is
easy to check that this procedure always produces a complex torus
with a flat K\"ahler metric and a B-field which is mirror to the
original one. This recipe is a special case of the physical notion
of T-duality~\cite{Polchtext}.

We will soon see how these results can be used to test the
Homological Mirror Conjecture for flat tori. First, let us recall
the formulation of the conjecture.

A physicist's Calabi-Yau $(X,G,\B)$ is both a complex manifold and
a symplectic manifold (the symplectic form being the K\"ahler form
$\omega=GI$). We can associate to each such manifold a pair of
triangulated categories: the bounded derived category of coherent
sheaves $\db{X}$ and the derived Fukaya category $\bD\cF_0(X).$
The former depends only on the complex structure of $X,$ while the
latter depends only on its symplectic structure. The Homological
Mirror Conjecture (HMC) asserts that if two algebraic Calabi-Yaus
$(X,G,\B)$ and $(X',G',\B')$ are mirror to each other, then
$\db{X}$ is equivalent to $\bD\cF_0(X'),$ and vice versa $\db{X'}$
is equivalent to $\bD\cF_0(X).$

Next, we need to recall the definitions of these two categories.
We begin with the Fukaya category $\cF(X).$ An object of the
Fukaya category is a triple$(Y,E,\nabla_E)$ where $Y$ is a graded
Lagrangian submanifold of $X,$ $E$ is a complex vector bundle on
$Y$ with a Hermitian metric, and $\nabla_E$ is a flat unitary
connection of $E.$ The only term which needs to be explained here
is ``graded Lagrangian submanifold.'' This notion was introduced
by J.~Leray in 1968 in the case of Lagrangian submanifolds in a
symplectic affine space, and generalized by Kontsevich to the
Calabi-Yau case in~\cite{Konts}. Let us first recall the
definition of the Maslov class of a Lagrangian submanifold. Let us
choose a holomorphic section $\Omega$ of the canonical line bundle
(which is trivial for Calabi-Yaus). Restricting it to $Y,$ we
obtain a nowhere vanishing $n$-form. On the other hand, we also
have a volume form $vol$ on $Y,$ which comes from the K\"ahler
metric on $X.$ This is also a nowhere vanishing $n$-form on $Y,$
and therefore $\Omega\vert_Y=h\cdot vol,$ where $h$ is a nowhere
vanishing complex function on $Y.$ The function $h$ can be thought
of as an element of $H^0(Y, {\mathcal C}^{\infty*}_Y),$ where
${\mathcal C}^{\infty*}_Y$ is the sheaf of $\CC^*$-valued
infinitely smooth functions on $Y.$ The standard exponential exact
sequence gives a homomorphism from $H^0(Y, {\mathcal
C}^{\infty*}_Y)$ to $H^1(Y,\ZZ),$ and the Maslov class of $Y$ is
defined as the image of $h$ under this homomorphism. (Explicitly,
the Cech cocycle representing the Maslov class is constructed as
follows: choose a good cover of $Y,$ take the logarithm of $h$ on
each set of the cover, divide by $2\pi i,$ and compare the results
on double overlaps). Although the definition of the Maslov class
seems to depend both on the complex and symplectic structures on
$X,$ in fact it is independent of the choice of complex structure.
Note also that if the Maslov class vanishes, the logarithm of $h$
exists as a function, and is unique up to addition of $2\pi i m,
m\in\ZZ.$ A graded Lagrangian submanifold $Y$ is a Lagrangian
submanifold in $X$ with a vanishing Maslov class and a choice of
the branch of $\log h.$

Morphisms are defined as follows. Suppose we are given two objects
$(Y_1,E_1,\nabla_1)$ and $(Y_2,E_2,\nabla_2).$ We will assume that
$Y_1$ and $Y_2$ intersect transversally at a finite number of
points; if this is not the case, we should deform one of the
objects by flowing along a Hamiltonian vector field, until the
transversality condition is satisfied. Let $\{e_i,i\in I\}$ be the
set of intersection points of $Y_1$ and $Y_2.$ Now we consider the
Floer complex. As a vector space, it is a direct sum of vector
spaces
$$
V_i=\Hom(E_1(e_i),E_2(e_i)),\quad i\in I.
$$
The grading is defined as follows. At any point $p\in Y$ the space
$T_pY$ defines a point $q$ in the Grassmannian of Lagrangian
subspaces in the space $T_pX.$ Let us denote by
$\widetilde{Lag}_p$ the universal cover of the Lagrangian
Grassmannian of the space $T_pX.$ On a Calabi-Yau variety $X,$
this spaces fit into a fiber bundle over $X$ denoted by
$\widetilde{Lag}$~\cite{Konts}. Grading of $Y$ provides a
canonical lift of $q$ to $\widetilde{Lag}_p$ for all $p$; this
lifts assemble into a section of the restriction of
$\widetilde{Lag}$ to $Y$~\cite{Konts}. Thus for each intersection
point $e_i$ we have a pair of points $q_1,q_2\in
\widetilde{Lag}_{e_i}.$ The grade of the component of the Floer
complex corresponding to $e_i$ is the Maslov index of $q_1,q_2$
(see ~\cite{Maslovindex} for a definition of the Maslov index.)
Finally, we need to define the differential. Let $e_i$ and $e_j$
be a pair of points whose grade differs by one. The component of
the Floer differential which maps $V_i$ to $V_j$ is defined by
counting holomorphic disks in $X$ with two marked points on the
boundary, so that the two marked points are $e_i$ and $e_j$ (the
Maslov index of $e_j$ is the Maslov index of $e_i$ plus one), and
the two intervals which make up the boundary of the disks are
mapped to $Y_1$ and $Y_2.$ Note that in order to compute the
differential one has to choose an (almost) complex structure $J$
on $X$ such that the form $\omega(\cdot,J\cdot)$ is a Hermitian
form on the tangent bundle of $X.$ For a precise definition of the
Floer differential, see~\cite{fooo}. The space of morphisms in the
Fukaya category is defined to be the Floer complex. The
composition of morphisms can be defined using holomorphic disks
with three marked points and boundaries lying on three Lagrangian
submanifolds.

The definition sketched above is only approximate. First, in order
to define the Floer differential one has to fix a relative spin
structure on $Y$~\cite{fooo}. Second, the Floer differential does
not square to zero in general, so the Floer ``complex'' is not
really a complex. A related difficulty is that the composition of
morphisms is associative only up to homotopy, which depends on
certain ternary product of morphisms. Actually, there is an
infinite sequence of higher products in the Fukaya category, which
are believed to satisfy the identities of an $A_\infty$ category
(see~\cite{Keller} for a review of $A_\infty$ categories). It is
also believed that changing the almost complex structure $J$ gives
an equivalent $A_\infty$ category, so that the equivalence class
of the Fukaya category is a symplectic invariant. For a detailed
discussion of these issues see~\cite{fooo}.

If the B-field is non-zero, one has to modify the definition of
the Fukaya category as follows. Objects are triples
$(Y,E,\nabla_E),$ where $Y$ is a graded Lagrangian submanifold,
$E$ is a vector bundle on $Y$ with a Hermitian metric, and
$\nabla_E$ is a Hermitian connection on $E$ such that its
curvature satisfies
$$
\nabla_E^2=2\pi i B\vert_Y.
$$
Thus the connection is projectively flat rather than flat.

Morphisms are modified in the following way: all occurrences of the
symplectic form $\omega$ in the definition of the Floer complex
and higher products are replaced with $\omega+iB.$ The modified
Fukaya category of a symplectic manifold $X$ with a B-field $B$
will be denoted $\cF(X,B).$

The Fukaya category $\cF(X,B)$ is not a true category, but an
$A_\infty$ category with a translation functor. The set of
morphisms between two objects in an $A_{\infty}$\!--category is a
differential graded vector space. To any $A_{\infty}$\!--category
one can associate a true category which has the same objects but
the space of morphisms between two objects is the $0$\!--th
cohomology group of the morphisms in the $A_{\infty}$\!--category.
Applying this construction to $\cF(X,B),$ we obtain a true
category $\cF_0(X,B)$ which is also called the Fukaya category.
Considering twisted complexes over $\cF(X,B)$
M.Kontsevich~\cite{Konts} also constructs a certain triangulated
category $\bD\cF_0(X,B).$
 We will call it the derived Fukaya category. The
category $\bD\cF_0(X,B)$ contains $\cF_0(X,B)$ as a full
subcategory.

In the next lecture we will discuss the derived category of
coherent sheaves, and use its properties to test the Homological
Mirror Conjecture.

\section{Derived Categories of Coherent Sheaves and a Test of the Homological Mirror Conjecture}

Let $X$ be a complex algebraic variety (or a complex manifold).
Denote by $\O_{X}$ the sheaf of regular functions (or the sheaf of
holomorphic functions). Recall that a coherent sheaf is a sheaf of
$\O_X$\!--modules that locally can be represented as a cokernel of
a morphism of algebraic (holomorphic) vector bundles. Coherent
sheaves form an abelian category which will be denoted by
$\coh(X).$

Next we recall the definition of a derived category and describe
some properties of derived categories of coherent sheaves on
smooth projective varieties. There is a lot of texts where
introductions to the theory of derived and triangulated categories
are given, we can recommend \cite{Ve, Ha, GM, KaSh, Ke}.

Let $\A$ be an abelian category. We denote by $\bC^b(\A)$ the
category of bounded differential complexes
$$
M^{\bdot}=(0\lto\cdots\lto M^p\stackrel{d^p}{\lto}
M^{p+1}\lto\cdots\lto 0), \quad M^p\in\A, \quad p\in\ZZ, \quad
d^2=0.
$$
Morphisms $f:M^{\bdot} \lto N^{\bdot}$ between complexes are sets
of morphisms $f^{p}:M^{p} \lto N^{p}$ in the category $\A$ which
commute with the differentials, i.e.
$$
d_N f^p - f^{p+1} d_M=0 \qquad\text{for all}\quad  p.
$$
A morphism of complexes $f: M^{\bdot}\lto N^{\bdot}$ is
null-homotopic if $f^p=d_N h^p + h^{p+1} d_M$ for all
 $p\in\ZZ$
for some family of morphisms $h^{p}: M^{p+1}\lto N^p.$ We define
the homotopy category $\bH^b(\A)$ as a category which has the same
objects as $\bC^b(\A),$ whereas morphisms in $\bH^b(\A)$ are the
equivalence classes of morphisms of complexes $\overline{f}$
modulo  the null-homotopic morphisms.

For any complex $N^{\bdot}$ and for any $p\in \ZZ$ we define the
cohomology $H^p(N^{\bdot})\in \A$ as the quotient $\Ker d^{p}/\Im
d^{p-1}.$ Hence for any  $p$ there is a functor
 $H^p : \bC(\A)\lto \A,$ which assigns to the complex $N^{\bdot}$
the cohomology $H^p(N^{\bdot})\in \A.$ We define a
quasi-isomorphism to be a morphism of complexes $s: N^{\bdot}\to
M^{\bdot}$ such that the induced morphisms $H^p s:
H^p(N^{\bdot})\to H^p(M^{\bdot})$ are invertible for all
$p\in\ZZ.$ We denote by $\Sigma$ the class of all
quasi-isomorphisms. This class of morphisms enjoys good properties
which are similar to the Ore conditions in the localization theory
of rings.

The bounded  derived category $\bD^b(\A)$ is defined as the
localization of $\bH^b(\A)$ with respect to the class $\Sigma$ of
all quasi-isomorphisms. This means that the derived category has
the same objects as the homotopy category $\bH^b(\A),$ and that
morphisms in the derived category from $N^{\bdot}$ to $M^{\bdot}$
are given by left fractions $s^{-1}\circ f,$ i.e. equivalence
classes of diagram
$$
N^{\bdot}\stackrel{f}{\lto}M^{'\bdot}\stackrel{s}{\longleftarrow}
M^{\bdot}, \qquad s\in \Sigma,
$$
where pairs $(f,s)$ and  $(g,t)$ are considered equivalent iff
there is a commutative diagram in $\bH^b(\A)$
$$
\xymatrix{
&M^{'\bdot}\ar[d]&\\
N^{\bdot}\ar[ur]^f \ar[r]^{h} \ar[dr]_{g}& M^{'''\bdot} &
M^{\bdot}\ar[ul]_s \ar[l]_{r}\ar[dl]^t\\
&M^{''\bdot}\ar[u]}
$$
such that $r\in\Sigma.$ Composition of morphisms $(f,s)$ and
$(g,t)$ is a morphism
 $(g'f, s't)$ which is defined using the commutative diagram:
$$
\xymatrix{
&&K^{''\bdot}&&\\
&M^{'\bdot}\ar@{-->}[ur]^{g'}&&K^{'\bdot}\ar@{-->}[ul]_{s'}&\\
N^{\bdot}\ar[ur]^{f}&&M^{\bdot}\ar[ul]_{s}\ar[ur]^g &&
K^{\bdot}\ar[ul]_t}.
$$
Such a diagram always exists, and one can check that the
composition law is associative.

We have a canonical functor $\bH^b(\A)\lto \bD^b(\A)$ sending a
morphism $f: N^{\bdot}\to M^{\bdot}$ to the pair $(f, {\rm
id}_{M}).$ This functor makes all quasi-isomorphisms invertible
and is universal among functors with this property. The abelian
category $\A$ can be considered as a full subcategory of
$\bD^{b}(\A)$ identifying an object $A\in \A$ with the complex
$\cdots\to 0\to A\to 0\to \cdots$ having $A$ in degree 0. If
$N^{\bdot}$ is an arbitrary complex, we denote by $N^{\bdot}[1]$
the complex with components $N^{\bdot}[1]^p=N^{p+1}$ and the
differential $d_{N[1]}= -d_N.$ This correspondence gives a functor
on the derived category $\bD^b(\A)$ which is an autoequivalence
and is called the translation functor.

Any derived category $\bD^b(\A)$ has a structure of a triangulated
category \cite{Ve}. This means  that the following data are
specified:
\begin{list}{\alph{tmp})}%
{\usecounter{tmp}} \item a translation functor $[1]:
\bD^b(\A)\lto\bD^b(\A)$ which is an additive autoequivalence,
\item a class of distinguished (or exact) triangles
$$
X\stackrel{u}{\lto}Y\stackrel{v}{\lto}Z\stackrel{w}{\lto}X[1]
$$
\end{list}
that satisfies a certain set of axioms (for details see~\cite{Ve,
GM, KaSh, Ke}).


To define a triangulated structure on the derived category
$\bD^b(\A)$ we introduce the notion of a standard triangle as a
sequence
$$
N\stackrel{Qi}{\lto} M\stackrel{Qp}{\lto} K\stackrel{\partial
\varepsilon}{\lto} N[1],
$$
where $Q: \bC^b(\A)\lto \bD^b(\A)$ is the canonical functor,
$$
0\lto N\stackrel{i}{\lto} M\stackrel{p}{\lto} K\lto 0
$$
is a short exact sequence of complexes, and $\partial \varepsilon$
is a certain morphism in $\bD^b(\A).$ The morphism $\partial
\varepsilon$ is the fraction $s^{-1} \circ j,$ where $j$ is the
inclusion of the subcomplex $K$ into the complex $C(p)$ with
components $K^n\oplus M^{n+1}$ and differential
$$
d_{C(p)}=
\begin{pmatrix}
d_K& p\\
0& -d_M
\end{pmatrix},
$$
and the quasi-isomorphism $s: N[1]\lto C(p)$ is the morphism $(0,
i).$

A {\it distinguished triangle} in $\bD^b(\A)$ is a sequence in
$\bD^b(\A)$ of the form
$$
X\stackrel{u}{\lto}Y\stackrel{v}{\lto}Z\stackrel{w}{\lto}X[1]
$$
which is isomorphic to a standard triangle.

Let $\A$ and $\cB$ be two abelian categories and  $F:\A\lto \cB$
an additive functor  which is left (resp. right) exact. The
functor $F$ induces a functor between the categories of
differential complexes and a functor $\bar{F}: \bH^b(\A)\to
\bH^b(\cB)$ obtained by applying $F$ componentwise. If $F$ is not
exact it does not transform quasi-isomorphisms into
quasi-isomorphisms. Nevertheless, often   we can define its right
(resp. left) derived functor ${\mathbf R} F$ (resp. ${\mathbf L}
F$) between the corresponding derived categories. The derived
functor ${\mathbf R} F$ (resp. ${\mathbf L} F$) will be  {\it
exact} functor between triangulated categories in the following
sense: it commutes with the translation functors and  takes every
distinguished triangle to a distinguished triangle. We will not
give here the definition of the derived functors, but the idea is
to apply the functor $F$ componentwise to well-selected
representatives of classes of quasi-isomorphic complexes
(see~\cite{Ve, Ha, GM, KaSh, Ke}).

For example, let us consider the derived categories of coherent
sheaves on smooth projective ( or proper algebraic) varieties. We
denote $\bD^b(\coh(X))$ by ${\db{X}}.$ Every morphism $f: X\to Y$
induces the inverse image functor $f^*: \coh(Y)\lto \coh(X).$ This
functor is right exact and has the left derived functor ${\mathbf
L}f^{*}: \db{Y}\lto \db{X}.$ To define it we have to replace a
complex on $Y$ by a quasi-isomorphic complex of locally free
sheaves, and apply the functor $f^*$ componentwise to this locally
free complex. Similarly, for any complex ${\F}\in {\db{X}}$ we can
define  an exact functor $\stackrel{\mathbf L}{\ot}\F :
{\db{X}}\lto{\db{X}}$ replacing $\F$ by a quasi-isomorphic locally
free complex.

The morphism $f: X\lto Y$ induces also the direct image functor
$f_*: \coh(X)\lto \coh(Y)$ which is left exact and has the right
derived functor ${\mathbf R}f_* : \db{X}\lto \db{Y}.$ To define it
we have to include the category of coherent sheaves into the
category of quasi-coherent sheaves and replace a complex by a
quasi-isomorphic complex of injectives. After that we can apply
the functor $f_*$ componentwise to the complex of injectives. The
functor ${\mathbf R} f_*$ is right adjoint to ${\mathbf L}f^*.$
This means that there is a functorial isomorphisms
$$
\Hom(A, {\mathbf R}f_* B)\cong \Hom({\mathbf L}f^* A, B).
$$
for all $A,B.$ This property can also be regarded as a definition
of the functor ${\mathbf R}f_*.$

Using these functors one can introduce a larger class of exact
functors. Let $X$ and $Y$ be  smooth projective (or proper
algebraic) varieties. Consider the projections
$$
\begin{array}{ccccc}
X&\stackrel{p}{\longleftarrow}&X\times Y&\stackrel{q}{\lto}&Y.\\
\end{array}
$$
Every object ${\E}\in {\db{X\times Y}}$ defines an exact functor
$\Phi_{\E}:{\db{X}}\lto{\db{Y}}$ by the following formula:
\begin{equation}\label{dfun}
\Phi_{\E}(\cdot):={\bf R}q_{*}({\E}\stackrel{\bf L}{\otimes}
p^*(\cdot)).
\end{equation}
Obviously, the same object defines another functor
$\Psi_{\E}:{\db{Y}}\lto{\db{X}}$ by a similar formula
$$
\Psi_{\E}(\cdot):={\bf R}p_{*}({\E}\stackrel{\bf L}{\otimes}
q^*(\cdot)).
$$
%

Thus there is a reasonably large class of exact functors between
bounded derived categories of smooth projective varieties that
consists of functors having the form $\Phi_{\E}$ for some complex
$\E$ on the product variety. This class is closed under
composition of functors. Indeed, let $X, Y, Z$ be three smooth
projective varieties and let
$$
\Phi_{I}:{\db{X}}\lto{\db{Y}},\qquad \Phi_{J}:{\db{Y}}\lto{\db{Z}}
$$
be two functors, where $I$ and $J$ are objects of ${\db{X\ts Y}}$
and ${\db{Y\ts Z}}$ respectively. Denote by $p_{XY}, p_{YZ},
p_{XZ}$ the projections of $X\ts Y\ts Z$ on the corresponding pair
of factors. The composition $\Phi_{J}\circ\Phi_{I}$ is isomorphic
to $\Phi_{K},$ where $K\in{\db{X\ts Z}}$ is given by the formula
$$
K\cong {\bf R}p_{XZ*}(p_{YZ}^*(J)\ot p_{XY}^*(I)).
$$


Presumably, the  class of exact functors described above
encompasses all exact functors between bounded derived categories
of coherent sheaves on smooth projective varieties. We do not know
if it is true or not. However it is definitely true for exact
equivalences.

\begin{theorem}{\rm (\cite{Or1},  also \cite{Orlovnew})} Let $X$ and $Y$
be smooth projective varieties. Suppose $F:
{\db{X}}\stackrel{\sim}{\lto}{\db{Y}}$ is an exact equivalence.
Then there exists a unique (up to an isomorphism) object
$\E\in{\db{X\ts Y}}$ such that the functor $F$ is isomorphic to
the functor $\Phi_{\E}.$
\end{theorem}


Now we consider the bounded derived categories of coherent sheaves
on abelian varieties. There are examples of different abelian
varieties which have equivalent derived categories of coherent
sheaves. Moreover, one can completely describe classes of abelian
varieties with equivalent derived categories of coherent sheaves.

Let $A$ be an abelian variety of dimension $n$ over $\C.$
 This means that $A$ is a complex torus $(U/\Gamma, I)$ which is
 algebraic, i.e. it has an embedding to the projective space.
Let $\wh{A}$ be the dual abelian variety, i.e. the dual torus
$(U^*/\Gamma^*, -I^t).$ It is canonically isomorphic to Picard
group $\Pic^0(A).$ It is well-known that there is a unique line
bundle $P$ on the product $A\ts\wh{A}$ such that for any point
$\alpha\in \wh{A}$ the restriction $P_{\alpha}$ on $A\ts
\{\alpha\}$ represents an element of $\Pic^0 A$ corresponding to
$\alpha,$ and, in addition, the restriction $P\Big|_{\{0\}\ts
\wh{A}}$ is trivial. Such $P$ is called the Poincare line bundle.

 The Poincare line bundle gives an example of an exact
equivalence between derived categories of coherent sheaves of two
non-isomorphic varieties. Let us consider the projections
$$
\begin{array}{ccccc}
A&\stackrel{p}{\longleftarrow}&A\times \wh{A}&\stackrel{q}{\lto}&\wh{A}\\
\end{array}
$$
and the functor $\Phi_{P}:{\db{A}}\lto{\db{\wh{A}}},$ defined as
in (\ref{dfun}), i.e. $\Phi_{P}(\cdot)={\bf R}q_{*}(P\ot
p^{*}(\cdot)).$ It  was proved by Mukai \cite{Mu1}  that the
functor $\Phi_{P}:{\db{A}}\lto{\db{\wh{A}}}$ is an exact
equivalence, and there is an isomorphism of functors:
$$
\Psi_{P}\circ\Phi_{P}\cong (-1_A)^*[n],
$$
where $(-1_A)$ is the inverse map on the group $A.$

Now let $A_1$ and $A_2$ be two abelian varieties of the same
dimension. We denote by $\Gamma_{A_1}$ and $\Gamma_{A_2}$ the
first homology lattices $H_1(A_1, \ZZ)$ and  $H_1(A_2,\ZZ).$ Every
map $f: A_1\lto A_2$ of  abelian varieties induces a map $\bar{f}:
\Gamma_{A_1}\lto \Gamma_{A_2}$ of the first homology groups.

For any abelian variety $A$ the first homology lattice of the
variety $A\times \wh{A}$ coincides with $\Gamma_A\oplus
\Gamma_A^*$ and hence it has the canonical symmetric bilinear form
$q_A$ defined by Equation~(\ref{bf}). Consider an isomorphism $f:
A_1\times \wh{A}_1\stackrel{\sim}{\lto} A_2\times \wh{A}_2$ of
abelian varieties. We call such map {\it isometric} if the
isomorphism $\bar{f}: \Gamma_{A_1} \oplus
\Gamma_{A_1}^*\stackrel{\sim}{\lto}\Gamma_{A_2} \oplus
\Gamma_{A_2}^*$ identifies the forms $q_{A_1}$ and $q_{A_2}.$

Now we can formulate a criterion for two abelian varieties to have
equivalent derived categories of coherent sheaves.

\begin{theorem}{\rm (\cite{Or2})}\label{eqcat} Let $A_1$ and $A_2$ be abelian
varieties. Then the derived categories $\db{A_1}$ and $\db{A_2}$
are equivalent as triangulated categories if and only if there
exists an isometric isomorphism
$$f :A_1\times \wh{A}_1\stackrel{\sim }{\lto}A_2\times \wh{A}_2,$$
i.e $\bar{f}$  identifies the forms $q_{A_1}$ and $q_{A_2}$ on
$\Gamma_{\vphantom{\wh{A}}A_1}\oplus \Gamma_{\wh{A}_1}$ and
$\Gamma_{\vphantom{\wh{A}}A_2}\oplus \Gamma_{\wh{A}_2}.$
\end{theorem}


Using Theorems~\ref{ntwoiso} and ~\ref{eqcat} we can now make a
check of the Homological Mirror Conjecture for tori. Suppose the
tori $(T_1,I_1,G_1,\B_1)$ and $(T_2,I_2,G_2,\B_2)$ are both mirror
to $(T',I',G',\B').$ Then $\SCFT(\Gamma_1,I_1,G_1,\B_1)$ is
isomorphic to $\SCFT(\Gamma_2,I_2,G_2,\B_2),$ and by
Theorem~\ref{ntwoiso} there is an isomorphism of lattices
$\Gamma_1\op\Gamma_1^*$ and $\Gamma_2\op\Gamma_2^*$ which
intertwines $q_1$ and $q_2,$ $\cI_1$ and $\cI_2,$ and $\cJ_1$ and
$\cJ_2.$

On the other hand, if we now assume that both complex tori
$(T_1,I_1)$ and $(T_2,I_2)$ are algebraic, then HMC implies that
$\db{(T_1,I_1)}$ is equivalent to $\db{(T_2,I_2)}.$ The criterion
for this equivalence is given in Theorem \ref{eqcat}: it requires
the existence of an isomorphism of $\Gamma_1\op\Gamma_1^*$ and
$\Gamma_2\op\Gamma_2^*$ which intertwines $q_1$ and $q_2,$ and
$\tilde{\cI_1}$ and $\tilde{\cI_2}.$ Clearly, since $\cI\neq
\tilde{\cI}$ in general, we get two conditions that   contradict
to each other. However, since $\cI$ coincides with $\tilde{\cI}$
under condition $\B^{(0,2)}=0$ we obtain the following result.
\begin{cor}\label{oneone}
If $\SCFT(\Gamma_1,I_1,G_1,\B_1)$ is isomorphic to
$\SCFT(\Gamma_2,I_2,G_2,\B_2),$ both $(T_1,I_1)$ and $(T_2,I_2)$
are algebraic, and both $\B_1$ and $\B_2$ are of type $(1,1),$
then $\db{(T_1,I_1)}$ is equivalent to $\db{(T_2,I_2)}.$
\end{cor}
Let $(T_1,I_1,G_1,B_1)$ be a complex torus equipped with a flat
K\"ahler metric and a B-field of type $(1,1)$ and let $(T_2,I_2)$
be another complex torus. Suppose there exists an isomorphism of
lattices $g:\Gamma_1\op\Gamma_1^*\ra \Gamma_2\op\Gamma_2^*$
mapping $q_1$ to $q_2$ and $\tilde{\cI}_1$ to $\tilde{\cI}_2.$ We
can prove that in this case there exists a K\"ahler metric $G_2$
and a B-field $B_2$ of type $(1,1)$ on $T_2$ such that
$\SCFT(\Gamma_1,I_1,G_1,\B_1)$ is isomorphic to
$\SCFT(\Gamma_2,I_2,G_2,\B_2)$ as an $N=2$ superconformal field
theory.

Combining this with Theorem~\ref{ntwoiso} and the criterion for
the equivalence of $\db{(T_1,I_1)}$ and $\db{(T_2,I_2)},$ we
obtain a result converse to Corollary~\ref{oneone}.
\begin{cor}
Let $(T_1,I_1,G_1,B_1)$ be an algebraic torus equipped with a flat
K\"ahler metric and a B-field of type $(1,1).$ Let $(T_2,I_2)$ be
another algebraic torus. Suppose $\db{(T_1,I_1)}$ is equivalent to
$\db{(T_2,I_2)}.$ Then on $T_2$ there exists a K\"ahler metric
$G_2$ and a B-field $B_2$ of type $(1,1)$ such that
$\SCFT(\Gamma_1,I_1,G_1,\B_1)$ is isomorphic to
$\SCFT(\Gamma_2,I_2,G_2,\B_2)$ as an $N=2$ superconformal field
theory.
\end{cor}

If $\dim_\CC T=1,$ then the B-field is automatically of type
$(1,1).$ Therefore the HMC passes the check in this special case.
Of course, this is what we expect, since the HMC is known to be
true for the elliptic curve~\cite{PZ}. On the other hand, for
$\dim_\CC T>1$ we seem to have a problem.

Not all is lost however, and a simple modification of the HMC
passes our check. The modification involves replacing $\db{(T,I)}$
with a derived category of $\beta(\B)$\!-twisted sheaves, where
$\beta(\B)$ is an element of $H^2((T, I), \O^*_T)$ depending on
the B-field $\B\in H^2(X,\RR).$

Let $X$ be an algebraic variety over $\CC,$ and let $\B\in
H^2(X,\RR).$ Consider the homomorphism $\beta: H^2(X, \RR)\to
H^2(X, \O^*_X)$ induced by the canonical map $\RR\lto \O^*_X$ from
the following commutative diagram of sheaves:
$$
\begin{array}{ccccccccc}
0&\lto&\ZZ&\lto&\RR&\llongrightarrow&\RR/\ZZ&\lto&0\\
&&\big\|&&\da&&\da&&\\
0&\lto&\ZZ&\lto&\O_X&\stackrel{\exp(2\pi i\cdot)}
{\llongrightarrow}&\O^*_X&\lto&0
\end{array}
$$
Any element $a\in H^2(X, \O_X^*)$ gives us an $\O^*_X$ gerbe
$\X_a$ over $X.$ Consider the category of weight-1 coherent
sheaves  $\coh_{1}(\X_a)$ on the gerbe $\X_a.$ Now our
triangulated category can be defined as the derived category
$\db{\coh_{1}(\X_{\beta(\B)})}$ which will be denoted as $\db{X,
\B}.$ Recall that weight-1 coherent sheaves on the gerbe $\X_a$
can be described as twisted coherent sheaves on $X$ in the
following way. Choose an open cover $\lbrace U_i \rbrace_{i\in I}$
of $X$ such that the element $a\in H^2(X, \O_X^*)$ is represented
by a \v{C}ech 2-cocycle  $a_{ijk}\in \Gamma(U_{ijk}, \O^*_X)$
where $U_{ijk}=U_i\cap U_j \cap U_k.$ Now an $a$\!-twisted  sheaf
can be defined as a collection of coherent sheaves $\F_i$ on $U_i$
for all $i\in I$ together with isomorphisms $ \phi_{ji}: \F_i
|_{U_{ij}}\stackrel{\sim}{\lto} \F_j |_{U_{ij}} $ for all $i,j \in
I$ (s.t. $\phi_{ij}=\phi_{ji}^{-1}$) satisfying the twisted
cocycle condition $\phi_{ij}\phi_{jk}\phi_{ki}=a_{ijk}\id.$

When $\beta(\B)$ is a torsion element of $H^2(X, \O^*_X),$ the
abelian category of twisted sheaves is equivalent to the abelian
category of coherent sheaves of modules over
 an Azumaya algebra $\A_{\B}$ which corresponds to this
element. This implies that the corresponding derived categories
are also equivalent.

Let us remind the definition and basic facts about Azumaya
algebras. Let $\A$ be an $\O_X$\!--algebra which is coherent as a
sheaf $\O_X$\!--modules.   Recall that $\A$ is called an Azumaya
algebra if it is locally free as a sheaf of $\O_X$\!--modules, and
for any point $x\in X$ the restriction $\A(x):=\A\otimes_{\O_X}
\CC(x)$ is isomorphic to a matrix algebra $M_r(\CC).$ A trivial
Azumaya algebra is an algebra of the form $\Endo(E)$ where $E$ is
a vector bundle. Two Azumaya algebras $\A$ and $\A'$ are called
similar (or Morita equivalent) if there exist vector bundles $E$
and $E'$ such that
$$
\A\otimes_{\O_X}\Endo(E)\cong \A'\otimes_{\O_X}\Endo(E').
$$
Denote by $\coh(\A)$ the abelian category of sheaves of (right)
$\A$\!--modules which are coherent as sheaves of
$\O_X$\!--modules, and by $\db{X, \A}$ the bounded derived
category of $\coh(\A).$ It is easy to see that if the Azumaya
algebras $\A$ and $\A'$ are similar, then the categories
$\coh(\A)$ and $\coh(\A')$ are equivalent, and therefore the
derived categories $\db{\A}$ and $\db{\A'}$ are equivalent as
well.

Azumaya algebras modulo Morita equivalence form a group with
respect to tensor product. This group is called the Brauer group
of the variety $X$ and is denoted by $Br(X).$ There is a natural
map
$$
Br(X)\lto H^2(X, \O^*_X).
$$
This map is an embedding and its image is contained in the torsion
subgroup $H^2(X, \O_X^*)_{tors}.$ The latter group is denoted by
$Br'(X)$ and called the cohomological Brauer group of $X.$ The
well-known Grothendieck conjecture asserts that  the natural map
$Br(X)\lto Br'(X)$ is an isomorphism for smooth varieties. This
conjecture has been proved for abelian
varieties~\cite{GroAbelian}.

Suppose now that $\beta(\B)$ is a torsion element of $H^2((T, I),
\O^*_T),$ and consider an Azumaya algebra $\A_{\B}$ which
corresponds to this element. The derived category $\db{(T,
I),\A_{\B}}$ does not depend on the choice of $\A_{\B}$ because
all these algebras are Morita equivalent. It can be shown that the
derived category $\db{(T, I), \B}$ is equivalent to the derived
category $\db{(T, I), \A_{\B}}.$

A sufficient condition for the equivalence of $\db{(T_1, I_1),
\B_1}$ and $\db{(T_2, I_2), \B_2}$ for the case of algebraic tori
is provided by the following theorem~\cite{Po}.

\begin{theorem}{\rm (\cite{Po})}\label{Poli}
Let $(T_1,I_1)$ and $(T_2,I_2)$ be two algebraic tori. Let
$\B_1\in H^2(T_1,\RR)$ and $\B_2\in H^2(T_2,\RR),$ and suppose
$\beta$ maps both $\B_1$ and $\B_2$ to torsion elements. If there
exists an isomorphism of lattices $\Gamma_1\op\Gamma_1^*$ and
$\Gamma_2\op\Gamma_2^*$ which maps $q_1$ to $q_2,$ and $\cI_1$ to
$\cI_2,$ then $\db{(T_1,I_1),\B_1}$ is equivalent to
$\db{(T_2,I_2),\B_2}.$
\end{theorem}
It appears plausible that this is also a necessary condition for
$\db{(T_1,I_1),\B_1}$ to be equivalent to $\db{(T_2,I_2),\B_2}.$
Combining Theorem~\ref{Poli} with our Theorem~\ref{ntwoiso}, we
obtain the following result.
\begin{cor}\label{cortwo}
Suppose $\SCFT(\Gamma_1,I_1,G_1,\B_1)$ is isomorphic to
$\SCFT(\Gamma_2,I_2,G_2,\B_2),$ both $(T_1,I_1)$ and $(T_2,I_2)$
are algebraic, and both $\B_1$ and $\B_2$ are mapped by $\beta$ to
torsion elements. Then $\db{(T_1,I_1),\B_1}$ is equivalent to
$\db{(T_2,I_2),\B_2}.$
\end{cor}
This corollary suggests that we modify the HMC by replacing
$\db{X}$ with $\db{X,\B}.$ Recall that in the presence of a
B-field the definition of the Fukaya category is modified, and
that the modified category is denoted $\bD\cF_0(X,B).$ The
modified HMC asserts that if $(X,G,\B)$ is mirror to
$(X',G',\B'),$ then $\db{X,\B}$ is equivalent to
$\bD\cF_0(X',\B').$ Corollary~\ref{cortwo} shows that this
conjecture passes the check which the original HMC fails.

In the case of the elliptic curve the modified HMC is not very
different from the original one. Since $h^{0,2}=0$ in this case,
the complex side is unaffected by the B-field, while on the
symplectic side its only effect is to complexify the symplectic
form (replacing $\omega$ with $\omega+iB$). For true Calabi-Yaus
(the ones whose holonomy group is strictly $SU(n)$ and not some
subgroup) $h^{0,2}$ also vanishes, and the complex side is again
unmodified, but on the symplectic side the effects of the B-field
can be rather drastic. For example, flat connections on Lagrangian
submanifolds must be replaced with projectively flat ones, and
this has the tendency to reduce the number of A-branes. But for
complex tori of dimension higher than one the B-field has
important effects on both A-branes and B-branes.

\section{The Category of A-branes and the Fukaya Category}

In this lecture we will discuss topological D-branes of type A
(A-branes) on Calabi-Yau manifolds. As it was stated above, the
set of A-branes has the structure of an additive category,
and if $X$ is mirror to $X',$ then the category of A-branes on $X$
should be equivalent to the category of B-branes on $X',$ and vice
versa. There is a lot of evidence that the category of B-branes on
$X$ is equivalent to $\bD^b(X)$\footnote{In this lecture we will
assume that the B-field is trivial, for
simplicity.}~\cite{Douglas,Laz2,AL,Dia,KatzSharpe}. The
Homological Mirror Conjecture is essentially equivalent to the
statement that the category of A-branes on $X$ is equivalent to
the derived Fukaya category of $X.$ As we will see below, this is
not true for some $X,$ so the Homological Mirror Conjecture needs
to be modified.

In the case when $X$ is an elliptic curve, the Homological Mirror
Conjecture, with some relatively minor modifications, has been
proved by Polishchuk and Zaslow in ~\cite{PZ}. On the other hand,
in ~\cite{KO2} it was shown that in general the Fukaya category is
only a full sub-category of the category of A-branes. In the case
when $X$ is a torus of dimension higher than two with a constant
symplectic form, we have constructed in ~\cite{KO2} some examples
of A-branes which are represented by coisotropic, rather than
Lagrangian submanifolds. So far we do not have a proposal how to
define the category of A-branes mathematically. The goal of this
lecture is to explain the results of ~\cite{KO2} and discuss the
many unresolved issues.

To show that in general the category of A-branes on $X$ is
``bigger'' than $\bD\cF_0(X),$ we will exhibit a certain mirror
pair $X$ and $X'$ such that the group $K^0(\bD^b(X))$ is strictly
bigger than $K^0$ of the Fukaya category $X'.$ In fact, to
simplify life, we will tensor $K^0(\bD^b(X))$ with $\QQ$ and use
the Chern character to map the rational K-theory to the cohomology
of $X.$ In the case of the derived category of coherent sheaves,
the Chern character $\ch$ takes values in the intersection of
$H^*(X,\QQ)$ and $\op_{p=0}^n H^{p,p}(X),$ which are both
subgroups of $H^*(X,\CC).$ (The Hodge conjecture says that the
image of $\ch$ should coincide with this intersection.) In the
case of the Fukaya category, one can say the following. Mirror
symmetry induces an isomorphism of $H^*(X,\CC)$ and $H^*(X',\CC),$
therefore the analogue of the Chern character for the Fukaya
category should take values in some subgroup of $H^*(X',\CC).$ A
natural candidate for the Chern character of an object
$(Y,E,\nabla)$ of the Fukaya category is the Poincar\'{e} dual of
the corresponding Lagrangian submanifold (taken over $\QQ$). This
conjecture can also be physically motivated, and we will assume it
in what follows.

Let $E$ be an elliptic curve, $e$ be an arbitrary point of $E,$
and $\End_e(E)$ be the ring of endomorphisms of $E$ which preserve
$e.$ For a generic $E$ we have $\End_e(E)=\ZZ,$ but for certain
special $E$ $\End_e(E)$ is strictly larger than $\ZZ.$ Such
special $E$'s are called elliptic curves with complex
multiplication. It is not difficult to check that $E$ has complex
multiplication if and only if its Teichm\"uller parameter $\tau$
is a root of a quadratic polynomial with integral coefficients.
Let $E$ be an elliptic curve with complex multiplication. Consider
an abelian variety $X=E^n,$ $n\geq 2.$ One can show that for such
a variety the dimension of the image of the Chern character
$$
\ch: K^0(\bD^b(X))\ot\QQ \longrightarrow H^*(X,\QQ)
$$
coincides with the intersection
$$
H^*(X,\QQ)\bigcap \left(\bigoplus_{p=0}^n H^{p,p}(X)\right)
$$
and has the dimension equal to
$$
\dim_\QQ \Im(\ch)=\binom{2n}{n}.
$$
Further, $X$ is related by mirror symmetry to a symplectic torus
$X'$ of real dimension $2n.$ Cohomology classes Poincar\'{e}-dual
to Lagrangian submanifolds in $X'$ lie in the kernel of the map
\begin{equation}\label{Lefmap}
\begin{CD}
H^n(X',\RR) @>{\wedge\omega}>> H^{n+2}(X',\RR).
\end{CD}
\end{equation}
This map is an epimorphism, and therefore the dimension of the
kernel is equal to $\binom{2n}{n}-\binom{2n}{n+2}.$ Thus the image
of the Chern character map for the Fukaya category of $X'$ has
dimension less or equal to $\binom{2n}{n}-\binom{2n}{n+2}.$ On the
other hand, the mirror relation between $X$ and $X'$ induces an
isomorphism of their cohomology groups~\cite{GLO}. If we make a
reasonable assumption that this isomorphism is compatible with the
equivalence of the categories of B-branes on $X$ and A-branes on
$X',$ we infer that the derived Fukaya category $\bD\cF_0(X')$
cannot be equivalent to the category of A-branes on $X'.$ On the
other hand, we expect on the physical grounds that the former is a
full sub-category of the latter.

This leaves us with a question: if not all A-branes are Lagrangian
submanifolds, how can one describe them geometrically? On the
level of cohomology, if the Chern character of A-branes does not
take values in the kernel of the map ~(\ref{Lefmap}), where {\it
does} it take values? In the case of flat tori, we can answer the
second question. In this case we know that a mirror torus is
obtained by dualizing a Lagrangian sub-torus, and can infer how
the cohomology classes transform under this operation. The answer
is the following~\cite{GLO,NC}. Suppose the original torus is of
the form $X=A\times B,$ where $A$ and $B$ are Lagrangian real
sub-tori, and the mirror torus is $X'=\hat{A}\times B.$ Consider a
torus $Z=A\times\hat{A}\times B.$ It has two obvious projections
$\pi$ and $\pi'$ to $X$ and $X'.$ On $A\times{\hat{A}}$ we also
have the Poincare line bundle $P$ whose Chern character will be
denoted $\ch(P).$ Using an obvious projection from $Z$ to
$A\times\hat{A},$ we may regard $\ch(P)$ as a cohomology class on
$Z.$ Given a cohomology class $\alpha\in H^*(X,\QQ),$ let us
describe its image under mirror symmetry~\cite{GLO}. We pull
$\alpha$ back to $Z$ using $\pi:Z\ra X,$ tensor with $\ch(P),$ and
then push forward to $X'$ using $\pi':Z\ra X'.$ This gives a
cohomology class $\alpha'\in H^*(X',\QQ)$ which is mirror to
$\alpha.$ The requirement that $\alpha$ be in the intersection of
$H^*(X,\QQ)$ and $\oplus_p H^{p,p}(X)$ implies the following
condition on $\alpha'$:
$$
\iota_{\omega^{-1}}\alpha'-\omega\wedge\alpha'=0.
$$
Here $\iota_{\omega^{-1}}$ is the operator of interior
multiplication by the bi-vector $\omega^{-1}.$ Cohomology classes
dual to Lagrangian submanifolds satisfy this condition, but there
are other solutions as well. For example, it is easy to construct
some solutions of the form $\alpha'=e^{a},$ where $a\in
H^2(X,\ZZ).$ This suggests that there exist line bundles on $X'$
which can be regarded as A-branes. We will see below that this
guess is correct.

To make further progress in understanding A-branes, we need to
rely on physical arguments. As explained in Lecture 1, a classical
A-brane is a boundary condition for a sigma-model with preserves
$N=2$ super-Virasoro algebra. In ~\cite{KO2} we analyzed this
condition assuming that an A-brane is described by the following
geometric data: a submanifold $Y$ in $X,$ a Hermitian line bundle
$E$ on $Y,$ and a unitary connection $\nabla_E$ on $E.$ We showed
that in order for a triple $(Y,E,\nabla_E)$ to be a classical
A-brane, the following three conditions are necessary and
sufficient.

\begin{itemize}
\item[(i)] $Y$ must be a coisotropic submanifold of $X.$ This
means that the restriction of the symplectic form $\omega$ to $Y$
must have a constant rank, and its kernel is an integrable
distribution $\LY\subset TY.$ (An equivalent definition: $Y$ is
coisotropic if and only if for any point $p\in Y$ the
skew-orthogonal complement of $TY_p\subset TX_p$ is contained in
$TY_p.$) We will denote by $\NY$ the quotient bundle $TY/\LY.$ By
definition, the restriction of $\omega$ to $Y$ descends to a
symplectic form $\sigma$ on the vector bundle $\NY.$ \item[(ii)]
The curvature 2-form $F=(2\pi i)^{-1} \nabla_E^2,$ regarded as a
bundle map from to $TY$ to $TY^*,$ annihilates $\LY.$ (The factor
$(2\pi i)^{-1}$ is included to make $F$ a real 2-form with
integral periods). This implies that $F$ induces on $\NY$ a
skew-symmetric pairing which we will denote $f.$ \item[(iii)] The
forms $\sigma$ and $f,$ regarded as maps from $\NY$ to $\NY^*,$
satisfy $(\sigma^{-1}f)^2=-\id_{\NY}.$ This means that
$\J=\sigma^{-1}f$ is a complex structure on the bundle $\NY.$
\end{itemize}

Let us make some comments on these three conditions. The condition
(i) implies the existence of a foliation of $Y$ whose dimension is
equal to the codimension of $Y$ in $X.$ It is known as the
characteristic foliation of $Y.$ If the characteristic foliation
happens to be a fiber bundle $p:Y\ra Z$ with a smooth base $Z,$
then $\NY$ is simply the pull-back of $TZ$ to $Y,$ i.e.
$\NY=p^*TZ,$ and the form $\sigma$ is a pull-back of a symplectic
form on $Z.$ In general, $\NY$ is a foliated vector bundle over
the foliated manifold $Y,$ and the space of leaves $Z$ is not a
manifold, or even a Hausdorff topological space. It still makes
sense to talk about the sheaf of local sections of $\NY$ locally
constant along the leaves of the foliation. This sheaf plays the
role of the pull-back of the sheaf of sections of the generally
non-existent tangent bundle to $Z.$ In the same spirit, the 2-form
$\sigma$ should be interpreted as a symplectic form on $Z.$ One
can summarize the situation by saying that $\sigma$ is a
transverse symplectic structure on a foliated manifold $Y.$ The
bundle $\NY$ is usually called the normal bundle of the foliated
manifold $Y$ (not to be confused with the normal bundle of the
submanifold $Y$ itself).

The condition (ii) says that for any section $v$ of $\LY$ we have
$\iota_v F=0.$ Since $dF=0,$ this implies that the Lie derivative
of $F$ along such $v$ vanishes, i.e. $F$ is constant along the
leaves of the foliation. In the case when the characteristic
foliation is a fiber bundle with a smooth base $Z,$ this is
equivalent to saying that $f$ is a pull-back of a closed 2-form on
$Z.$ In general, one can say that $f$ is a transversely-closed
form on a foliated manifold $Y.$

The condition (iii) implies, first of all, that $f$ is
non-degenerate. Thus $f$ is a transverse symplectic structure on
$Y,$ just like $\sigma.$ Second, the condition (iii) says that the
ratio of the two transverse symplectic structures is a transverse
almost complex structure on the foliated manifold $Y.$ If the
characteristic foliation is a fiber bundle with a smooth base $Z,$
then $\J=\sigma^{-1}f$ is simply an almost complex structure on
$Z.$

An easy consequence of these conditions is that the dimension of
$Y$ must be $n+2k,$ where $n=\frac{1}{2}\dim_\RR X,$ and $k$ is a
non-negative integer. The number $k$ has the meaning of
``transverse complex dimension'' of $Y.$ If $k=0,$ then the
condition (i) says that $Y$ is a Lagrangian submanifold, and then
the condition (ii) forces $F$ to vanish. (The condition (iii) is
vacuous in this case). Another interesting special case occurs
when $Y=X$ (this is possible only if $n$ is even). In this case
the leaves of the characteristic foliation are simply points, the
conditions (i) and (ii) are trivially satisfied, and the condition
(iii) says that $\omega^{-1}F$ is an almost complex structure on
$X.$

A less obvious property is that the transverse almost complex
structure $\J$ is integrable~\cite{KO2}. This follows easily from
the well-known Gelfand-Dorfman theorem which plays an important
role in the theory of integrable systems~\cite{GD}. Thus $Y$ is a
transverse complex manifold. It is also easy to see that both $f$
and $\sigma$ have type $(0,2)+(2,0)$ with respect to $\J.$ In
fact, $f+i\sigma$ is a transverse holomorphic symplectic form on
the transverse complex manifold $Y.$

The somewhat mysterious condition (iii) can be rewritten in several equivalent forms.
For example, an equivalent set of
conditions is
$$
\wedge^r (f+i\sigma)\neq 0, \quad r< k, \quad \wedge^k (f+i\sigma)=0.
$$
Here $k$ is related to the dimension of $Y$ as above.

Our attempts to generalize the conditions (i)-(iii) to A-branes which carry vector bundles of rank higher than
one have been only partially successful. The first two conditions (i) and (ii) remain unchanged, but the
condition (iii) causes problems. Both physical and mathematical arguments indicate that
the correct generalization of (iii) looks as follows:
$$
(\sigma^{-1}f)^2=-{\rm id}_{E\ot\NY}.
$$
Here the ``transverse'' curvature 2-form $f$ is regarded as a
section of $\End(E)\ot\Lambda^2(\NY^*).$ This condition on $f$
does not lead to a transverse complex structure on $Y,$ and its
geometric significance is unclear. This leads to problems when one
tries to understand morphisms between such A-branes (see below).

So far our discussion of A-branes was classical. One the quantum
level $N=2$ super-Virasoro can be broken by anomalies. The absence
of such anomalies is an important additional condition on
A-branes. Let us focus our attention on the R-current $J$ whose
Fourier modes we denoted by $J_n$ and $\bJ_n$ previously. Its
conservation can be spoiled only by non-perturbative effects on
the world-sheet, i.e. by Riemann surfaces in $X$ whose boundaries
lie in $Y$ and which cannot be continuously deformed to a point.
In the case when $Y$ is Lagrangian, the conditions for the absence
of anomalies have been analyzed by K.~Hori in ~\cite{book}. The
result is that there are no anomalies if and only if the Maslov
class of $Y$ vanishes. This provides a physical interpretation of
the vanishing of the Maslov class for objects of the Fukaya
category.

For coisotropic $Y$ the condition for anomaly cancellation has
been obtained in ~\cite{Li}. Let $F$ be the curvature 2-form of
the line bundle on $Y,$ and let the dimension of $Y$ be $n+2k,$ as
before. Let $\Omega$ be a holomorphic trivialization of the
canonical line bundle on $X.$ One can show that the $n+2k$-form
$\Omega\vert_Y\wedge F^k$ is nowhere vanishing, and therefore we
have $\Omega\vert_Y\wedge F^k=h\cdot vol,$ where $vol$ is the
volume form, and $h$ is a smooth nowhere-vanishing complex-valued
function on $Y.$ We may regard $h$ as an element of $H^0(Y,
{\mathcal C}^{\infty *}_Y),$ where ${\mathcal C}^{\infty *}_Y$ is
the sheaf of smooth $\CC^*$-valued functions on $Y.$ Let
$\alpha_h\in H^1(Y,\ZZ)$ be the image of $h$ under the Bockstein
homomorphism. The anomaly of the R-current is absent if and only
if $\alpha_h=0.$

In the case when $X$ is a torus with a constant symplectic form, $Y$ is an affine sub-torus, and
the curvature 2-form $F$ is constant, one can quantize the sigma-model and verify directly  that the $N=2$
super-Virasoro algebra is preserved on the quantum level. This shows that non-Lagrangian A-branes exist
on the quantum level.

We hope that we have given convincing arguments that A-branes are not necessarily associated to
Lagrangian submanifolds, and that the Fukaya
category should be enriched with more general coisotropic A-branes.  Unfortunately, we do not have a
definite proposal for what should replace the
Fukaya category. In the remainder of this lecture we will describe a heuristic idea which could help to solve this problem.

In these lectures we have encountered two kinds of A-branes.
First, we have discussed objects of the Fukaya category, i.e.
triples $(Y,E,\nabla_E)$ where $Y$ is a graded Lagrangian
submanifold, $E$ is a vector bundle on $Y$ with a Hermitian
metric, and $\nabla_E$ is a flat unitary connection on $E.$
Second, we have triples $(Y,E,\nabla_E),$ where $E$ is a Hermitian
{\it line} bundle on $Y,$ $\nabla_E$ is a unitary connection on
$E,$ and the conditions (i)-(iii) are satisfied. For objects of
the Fukaya category we know in principle how to compute spaces of
morphisms and their compositions. Let us try to guess what the
recipe should be for objects of the second type.

We need to generalize the Floer complex to coisotropic A-branes.
To guess the right construction, it is useful to recall the
intuition which underlies the definition of the Floer complex.
Consider the space of smooth paths in $X,$ which we will denote
$PX.$ This space is infinite-dimensional, but let us treat it as
if it were a finite-dimensional manifold. We have a natural 1-form
$\alpha$ on $PX$ obtained by integrating the symplectic form
$\omega$ on $X$ along the path. More precisely, if $\gamma:I\ra X$
is a path, and $\beta$ is a tangent vector to $PX$ at point
$\gamma$ (i.e. a vector field along the image of $\gamma$), then
the value of $\alpha$ on $\beta$ is defined to be
$$
\int_I \omega(\stackrel{\cdot}\gamma(t),\beta(t)) dt.
$$
Note that the space $PX$ has two natural projections to $X$ which
we denote $\pi_1$ and $\pi_2.$ It is easy to see that
$d\alpha=\pi_2^*\omega-\pi_1^*\omega.$ Thus if we consider a
submanifold in $PX$ consisting of paths which begin and end on
fixed Lagrangian submanifolds in $X,$ then the restriction of
$\alpha$ to such a submanifold will be closed.

Let $Y_1$ and $Y_2$ be Lagrangian submanifolds in $X,$ and let
$PX(Y_1,Y_2)$ be the submanifold of $PX$ consisting of paths which
begin at $Y_1$ and end at $Y_2.$ Then the operator
$Q=d+2\pi\alpha$ on the space of differential forms on
$PX(Y_1,Y_2)$ satisfies $Q^2=0,$ and we may try to compute its
cohomology (in the finite-dimensional case this complex is called
the twisted de Rham complex.) Since $PX(Y_1,Y_2)$ is
infinite-dimensional, it is not easy to make sense of the twisted
de Rham cohomology. A.~Floer solved this problem by a formal
application of Morse theory. That is, if we formally apply the
Morse-Smale-Witten-Novikov theory to the computation of the
cohomology of $Q,$ we get precisely the Floer complex for the pair
$Y_1,Y_2.$

This construction ignores the bundles $E_1$ and $E_2,$ but it is
easy to take them into account. Consider the bundles $\pi_1^*
E_1^*$ and $\pi_2^* E_2$ on $PX(Y_1,Y_2).$ Both vector bundles
have natural unitary connections obtained by pulling back
$\nabla_{E_1^*}$ and $\nabla_{E_2}.$ We tensor them, and then add
the 1-form $2\pi\alpha$ to the connection on the tensor product.
Since $d\alpha=0,$ the resulting connection is still flat, but no
longer unitary. If we formally use Morse theory to compute the
cohomology of the resulting twisted de Rham complex on
$PX(Y_1,Y_2),$ we get the Floer complex for a pair of objects of
the Fukaya category $(Y_1,E_1,\nabla_{E_1})$ and
$(Y_2,E_2,\nabla_{E_2}).$

Now consider a pair of coisotropic A-branes, instead of a pair of
Lagrangian A-branes. We assume in addition that the bundles $E_1$
and $E_2$ are line bundles. By $PX(Y_1,Y_2)$ we still denote the
space of smooth paths in $X$ beginning at $Y_1$ and ending on
$Y_2.$ The first difficulty is that the restriction of $\alpha$ to
$PX(Y_1,Y_2)$ is not closed, so we cannot use it to define a
twisted de Rham complex. The second difficulty is that connections
$\pi_1^*(\nabla_{E_1^*})$ and $\pi_2^*(\nabla_{E_2})$ on the
bundles $\pi_1^* E_1^*$ and $\pi_2^* E_2$ are also not flat.
However, thanks to conditions (ii) and (iii) these two
difficulties ``cancel'' each other, as we will see in a moment. So
let us proceed as in the Lagrangian case: tensor the line bundles
$\pi_1^* E_1^*$ and $\pi_2^* E_2,$ and add $2\pi\alpha$ to the
connection. The resulting connection is not flat, but it has the
following interesting property. Note that since $Y_1$ and $Y_2$
are foliated manifolds, $PX(Y_1,Y_2)$ also has a natural
foliation. A leaf of this foliation consists of all smooth paths
in $X$ which begin on a fixed leaf in $Y_1$ and end on a fixed
leaf of $Y_2.$ The codimension of the foliation is finite and
equal to the sum of codimensions of characteristic foliations of
$Y_1$ and $Y_2.$ Further, since $Y_1$ and $Y_2$ have natural
transverse complex structures, the foliated manifold $PX(Y_1,Y_2)$
also has one. The connection on $\pi_1^* E_1^*\ot\pi_2^* E_2$ has
the following properties:

\begin{itemize}
\item[(A)] it is flat along the leaves of the foliation;
\item[(B)] its curvature $\pi_1^*(F_1+i\omega)-\pi_2^*(F_2+i\omega)$ has type $(2,0)$ in the transverse directions.
\end{itemize}

Thus it makes sense to consider a sheaf of sections of the line
bundle $\pi_1^* E_1^*\ot\pi_2^* E_2$ which are covariantly
constant along the leaves and holomorphic in the transverse
directions. It is natural to propose the cohomology of this sheaf
as the candidate for the space of morphisms between the A-branes
$(Y_1,E_1,\nabla_{E_1})$ and $(Y_2,E_2,\nabla_{E_2}).$ This
proposal can be formally justified by considering the path
integral quantization of the topologically twisted $\sigma$-model
on an interval with boundary conditions corresponding to the
A-branes $(Y_1,E_1,\nabla_{E_1})$ and $(Y_2,E_2,\nabla_{E_2}).$ To
get a rigorous definition of the spaces of morphisms, our formal
proposal must be properly interpreted. In the case of Lagrangian
A-branes, the above sheaf becomes the sheaf of covariantly
constant sections of a flat line bundle on $PX(Y_1,Y_2),$ and one
can interpret its cohomology using Morse-Smale-Witten-Novikov
theory. We do not know how to make sense of our formal proposal in
general.

The difficulty of generalizing the Floer complex to coisotropic A-branes suggests that perhaps the
geometric description of A-branes by means of submanifolds and vector bundles on them is not the right
way to proceed. Let us explain what we mean by this using an analogy from complex-analytic geometry.
There exists a general notion of a holomorphic vector bundle on a complex manifold, whose special
case is the notion of a holomorphic line bundle. One can study
line bundles in terms of their divisors, but this approach does not extend easily to higher rank bundles. Perhaps objects
of the Fukaya category, as well as coisotropic A-branes of rank one, are symplectic analogues of divisors, and in order to
make progress one has to find a symplectic analogue of the notion of a holomorphic vector bundle (or a coherent sheaf). This
analogy is strengthened by the fact that both divisors and geometric representations of A-branes by means of Lagrangian
or coisotropic submanifolds provide a highly
redundant description of objects in the respective categories: a line bundles does not change if one
adds to the divisor a principal divisor, while objects of the Fukaya category are unchanged by flows along Hamiltonian vector fields.
We believe that a proper definition of the category of A-branes will be very useful for understanding
Mirror Symmetry, and perhaps also for symplectic geometry as a whole.

\end{document}